\input ./style/arxiv-general.cfg
\documentclass[MSNbibl,number,citesort,seceqn,nameyear,rotating,dvips]{arxbj}
\makeatletter
   \@ifpackageloaded{graphicx}{}{\usepackage{graphicx}}
\makeatother

\volume{22}
\issue{2}
\pubyear{2016}
\firstpage{857}
\lastpage{900}
\doi{10.3150/14-BEJ678} 
\docsubty{FLA}

\makeatletter
\newcommand{\lleft}{\left}
\newcommand{\rright}{\right}
\newcommand{\rrVert}{\Vert}
\newcommand{\llVert}{\Vert}
\renewcommand{\mid}{|}
\newtheorem{lem}{Lemma}[section]
\newtheorem{teo}{Theorem}[section]
\newtheorem{cor}{Corollary}[section]
\def\real{\mathop{\mathbb R}} 
\def\nat{\mathop{\mathbb N}} 
\def\argmin{\mathop{\arg\min}}
\def\argmax{\mathop{\arg\max}}
\newcommand{\hatgnh}{\hat{g}_n}
\newcommand{\hathnh}{\hat{h}_n}
\newcommand{\hatfnh}{\hat{f}_n}
\newcommand{\hatmnh}{\hat{m}_n}
\newcommand{\cnxgb}{c_{n\bar{g}}}
\newcommand{\cnxmb}{c_{n\bar{m}}}
\newcommand{\calX}{\mathcal{X}}
\makeatother

\begin{document}
\begin{frontmatter}

\title{Consistency, efficiency and robustness of conditional disparity methods}
\runtitle{Conditional disparity methods}

\begin{aug}
\author[A]{\inits{G.}\fnms{Giles}~\snm{Hooker}\corref{}\ead[label=e1]{Giles.hooker@cornell.edu}}
\address[A]{Department of Biological Statistics and Computational Biology,
Cornell University, Ithaca, NY 14853-4201, USA. \printead{e1}}
\end{aug}

\received{\smonth{4} \syear{2014}}
\revised{\smonth{8} \syear{2014}}

%
\begin{abstract}
This paper considers extensions of minimum-disparity estimators to the
problem of estimating
parameters in a regression model that is conditionally specified; that
is where a parametric model describes the distribution of a response
$y$ conditional on covariates $x$ but does not specify the distribution
of $x$. We define these estimators by estimating a non-parametric
conditional density estimates and minimizing a disparity between this
estimate and the parametric model averaged over values of $x$. The
consistency and asymptotic normality of such estimators is demonstrated
for a broad class of models in which response and covariate vectors can
take both discrete and continuous values and incorportates a wide set
of choices for kernel-based conditional density estimation. It also
establishes the robustness of these estimators for a broad class of
disparities. As has been observed in
Tamura and Boos (\textit{J.~Amer. Statist. Assoc.} \textbf{81} (1986) 223--229), minimum
disparity estimators incorporating kernel density estimates of more
than one dimension can result in an asymptotic bias that is larger that
$n^{-1/2}$ and we characterize a similar bias in our results and show
that in specialized cases it can be eliminated by appropriately
centering the kernel density estimate. We also demonstrate empirically
that bootstrap methods can be employed to reduce this bias and to
provide robust confidence intervals. In order to demonstrate these
results, we establish a set of $L_1$-consistency results for
kernel-based estimates of centered conditional densities.
\end{abstract}

%
\end{frontmatter}

\section{Introduction}

Minimum disparity estimators (MDEs) are based on minimizing a measure
of distance between a non-parametric density estimate $\hatfnh(y)$ and
a parametric family of densities $\phi_{\theta}(y)$. Disparities can
be written in the general form Lindsay \cite{Lindsay94}:
\[
D(\hatfnh,\theta) = \int C \biggl( \frac{\hatfnh(y) - \phi
_{\theta
}(y)}{\phi_{\theta}(y)} \biggr)
\phi_{\theta}(y) \,\mathrm{d}\nu(y),
\]
where $C$ is a convex function with a minimum at 0 and $\nu$ is a
reference measure over the space of $y$. The minimum disparity
estimator is defined to be
\[
\hat{\theta}_n = \argmin_{\theta} D(\hatfnh,\theta).
\]
When $\hatfnh$ is a kernel density estimate based on univariate i.i.d.
data and $C(\delta)$ behaves appropriately at 0, these estimators can
be shown to be asymptotically normal and efficient in the sense of
having asymptotic variance given by the inverse of the Fisher
information. When $C$ behaves appropriately at $\infty$, they are also
robust to outliers. This was first observed in the case of Hellinger distance ($C(\delta) = [\sqrt{\delta+1}-1 ]^2$) by
Beran \cite{Beran77} and generalized to the broader class of
disparities in Lindsay \cite{Lindsay94} for discrete data and for
continuous data in Basu and Lindsay  \cite{BasuLindsay94} and Park and Basu \cite
{ParkBasu04}. The particular case of $C(\delta) = \mathrm{e}^{-\delta}$ was
studied in
Basu, Sahadeb and Vidyashankar
\cite{BSV97}; a choice that that is both robust to
outliers and to ``inliers'' -- regions where $\delta(\cdot) =
[\hatfnh(x\cdot) - \phi_{\theta}]/\phi_{\theta
}$ is near it's
negative limit of $-$1 and where Hellinger distance performs poorly.
Tamura and Boos \cite{TamuraBoos86} observed that when $\hatfnh(x\cdot)$ is a
multivariate kernel density estimate, the MDE has an asymptotic bias
that is larger than $n^{-1/2}$ and hence appears in the central limit
theorem for $\hat{\theta}_n$, potentially necessitating a bias correction.

Despite the potential for both robust and efficient estimation, minimum
disparity estimation has seen few extensions beyond i.i.d. data. Within
this context, the use of disparity methods to estimate parameters in
linear regression was treated in Pak and Basu \cite{PakBasu98} by
placing a disparity on the score equations and for discrete covariates
in Cheng and Vidyashankar \cite{ChenVidyashankar06}, but little attention has been
given to more general regression problems and we take a more direct
approach here. In this paper, we consider data
$(X_1,Y_1),(X_2,Y_2),\ldots$ for which we have a parameterized family
of densities $\phi_{\theta}(y\mid x)$ which describe the distribution of
$y$ conditional on the value of $x$. We construct a non-parametric
conditional density estimate $\breve{f}_n(y\mid x)$ based on kernel densities
and define two extensions of disparities:
\begin{eqnarray*}
D_n(\breve{f}_n,\theta) &=& \frac{1}{n} \sum
_{i=1}^n D\bigl(\breve{f}_n(
\cdot\mid X_i),\phi_{\theta}(\cdot\mid X_i)
\bigr),
\\
\tilde{D}_n(\breve{f}_n,\theta) &=& \int D\bigl(
\breve{f}_n(\cdot\mid x),\phi_{\theta}(\cdot\mid x)\bigr)
\hathnh(x)\,\mathrm{d}x,
\end{eqnarray*}
where $\hathnh(x)$ is a kernel density estimate of the density of $x$.
We show that the parameters minimizing these disparities are consistent
and asymptotically normal. Furthermore, when the data are generated
from a process that corresponds to some member of the parametric model,
the limiting variance is given by the information matrix. Our framework
is intentionally general and designed to cover a broad range of cases
in which both $Y_i$ and $X_i$ can be vector valued and incorporate a
mix of continuous- and discrete-valued components and are designed to
be as general as possible. We also consider various estimates of
$\breve{f}_n(y\mid x)$ in which some components of $y$ are centered by a
Nadaraya--Watson estimator based on some components of $x$. When the
parametric model is incorrect, these yield different bias and variance
expressions in our central limit theorem which we interpret and describe.

To achieve these results, we first demonstrate the $L_1$ consistency of
$\breve{f}_n(\cdot\mid x)$ which holds uniformly over $x$. We also
demonstrate the robustness of these estimators to outlying values in
$y$. The effectiveness of these techniques are then examined in
simulation and with real-world data.

We will introduce the specific distributional framework and assumptions
in the next subsection and our conditional density estimators in
Section~\ref{seckernels}. Because of the notational complexity
involved with working with both continuous and discrete random
variables as well as a division of the components of $x$, Section~\ref
{secnotation} will detail notational shorthand that will be used in
various places throughout the remainder of the paper. Section~\ref
{consistency1} will develop results on the $L_1$ consistency of
kernel-based conditional density estimators, Section~\ref
{secconsistentdisparity} will then apply these results to demonstrate
the consistency of minimum-disparity estimators in conditionally
specified models. We will demonstrate the asymptotic normality of these
estimators in Section~\ref{secCLT} and their robustness will be
examined in Section~\ref{secrobustness}. Computational details on
selecting bandwidths and using the bootstrap for bias correction and
inference are given in Section~\ref{secbootstrap}. Simulation results
and real data analysis are given in Sections~\ref{secsimulations} and
\ref{secdata}.

We have included proofs of our results in the text where they are
either enlightening or short, but have reserved many for a Supplemental Appendix (Hooker \cite{Hooker14}) and noted where these may be found.

\subsection{Framework and assumptions} \label{secassumptions}

Throughout the following, we assume a probability space $(\Omega
,\mathcal{F},P)$ from which
we observe i.i.d. random variables
$\{X_{n1}(\omega),X_{n2}(\omega
),Y_{n1}(\omega),Y_{n2}(\omega),n\geq1\}$
where we have separated discrete and continuous random variables so
that $X_{n1}(\omega) \in\real^{d_x}$,
$X_{n2}(\omega) \in S_x$, $Y_{n1}(\omega) \in\real^{d_y}$,
$Y_{n2}(\omega) \in
S_y$ for countable sets $S_x$ and $S_y$ with joint distribution
\[
g(x_1,x_2,y_1,y_2) =
P(X_2 = x_2,Y_2=y_2)P(X_1
\in\mathrm{d}x_1,Y_1 \in \mathrm{d}y_1\mid
X_2=x_2,Y_2=y_2)
\]
and define the marginal and conditional densities
%
\begin{eqnarray}
h(x_1,x_2) & =& \sum_{y_2 \in S_y}
\int g(x_1,x_2,y_1,y_2)
\,\mathrm{d}y_1,
\\
f(y_1,y_2\mid x_1,x_2) & =&
\frac{g(x_1,x_2,y_1,y_2)}{h(x_1,x_2)} \label{deff}
\end{eqnarray}
on the support of $(x_1,x_2)$.

An important aspect of this paper is to study an approach of centering
$y_1$ by a Nadaraya--Watson estimator before estimating $g$. We define
this generally, so that $y_1$ can be centered based on some components
$(X^{\bar{m}}_1,X^{\bar{m}}_2)$ of $(X_1,X_2)$ and a density for the
residuals can be estimated based on a different possibly-overlapping
set of components $(X^{\bar{g}}_1,X^{\bar{g}}_2)$. Formally, we
define $(x^{\bar{m}}_1,x^{\bar{m}}_2)$ and $(x^{\bar{g}}_1,x^{\bar
{g}}_2)$ with\vspace*{1pt} densities $h^{\bar{m}}(x^{\bar{m}}_1,x^{\bar{m}}_2)$
and $h^{\bar{g}}(x^{\bar{g}}_1,x^{\bar{g}}_2)$, respectively, where
$x^{\bar{m}}_1 \in\real^{d_{x\bar{m}}}$ and $x^{\bar{g}}_1 \in
\real
^{d_{x\bar{g}}}$ and $x^{\bar{m}}_2 \in S_{x^{\bar{m}}}$, $x^{\bar{g}}_2
\in S_{x^{\bar{g}}}$. We now define the possibly vector-valued
expectation of~$y_1$ conditional on $x^{\bar{m}}$:
\[
m\bigl(x^{\bar{m}}_1,x^{\bar{m}}_2\bigr) =
\sum_{y_2 \in S_y} \sum_{x^{\bar{g}}_2 \in S_{x^{\bar{g}}}} \int\!\!\!\int y_1 \frac{g(x_1,x_2,y_1,y_2)}{h^{\bar{m}}(x^{\bar{m}}_1,x^{\bar
{m}}_2)} \,\mathrm{d}y_1 \,\mathrm{d}x^{\bar{g}}_1
\]
along with the residuals
\[
\varepsilon= y_1 - m\bigl(x^{\bar{m}}_1,x^{\bar{m}}_2
\bigr)
\]
and define the joint density of these residuals, $y_2$, and $x^{\bar
{g}}$ by
\[
g^c\bigl(x^{\bar{g}}_1,x^{\bar{g}}_2,
\varepsilon,y_2\bigr) = \sum_{x^{\bar
{m}}_2 \in S_{x^{\bar{m}}}} \int g
\bigl(x_1,x_2,\varepsilon+m\bigl(x^{\bar
{m}}_1,x^{\bar{m}}_2
\bigr),y_2\bigr) \,\mathrm{d} x^{\bar{m}}_1
\]
and similarly write the conditional density
\[
f^c\bigl(\varepsilon,y_2\mid x^{\bar{g}}_1,x^{\bar{g}}_2
\bigr)= \frac{g^c(x^{\bar
{g}}_1,x^{\bar{g}}_2,\varepsilon,y_2) }{h^{\bar{g}}(x^{\bar
{g}}_1,x^{\bar{g}}_2)},
\]
where throughout this paper we will assume that the distribution of
$(y_1,y_2)$ is such that
\[
f(y_1,y_2\mid x_1,x_2) =
f^c\bigl(\varepsilon+m\bigl(x^{\bar{m}}_1,x^{\bar
{m}}_2
\bigr),y_2\mid x^{\bar{g}}_1,x^{\bar{g}}_2
\bigr)
\]
for some function $f^c(\varepsilon,y_2 \mid
x^{\bar{g}}_1,x^{\bar{g}}_2)$
that does not depend on those components of $X$ that are not also
components of $X^{\bar{g}}$.

A useful example to keep in mind is the conditionally heteroscedastic
linear regression model
\[
y_i = \bigl(x^{\bar{m}}_i \bigr)^T
\beta+ \sigma\bigl( \bigl(x^{\bar{g}}_i \bigr)^T
\gamma\bigr) \varepsilon
\]
for $\varepsilon\sim f(\cdot)$ in which the residual variance depends on
covariates $x^{\bar{g}}$ while the mean depends on $x^{\bar{m}}$ and
these may or may not be the same variables. However, our framework is
considerably more general than this model and includes all of ANOVA,
multiple regression, ANCOVA, multivariate regression, tabular data and
generalized linear models as well as allowing for more complex models
in which dependence is assumed between categorical and continuous
response variables.

To appreciate the generality of class of conditional density estimates,
we observe that this covers the case (\ref{deff}) by setting the
collection of variables in $(x^{\bar{m}}_1,x^{\bar{m}}_2)$ to be
empty and $(x^{\bar{g}}_1,x^{\bar{g}}_2) = (x_1,x_2)$; in this case
we understand $m(x^{\bar{m}}_1,x^{\bar{m}}_2) \equiv0$. It also
covers the ``homoscedastic'' in which there is no $y_2$ and we assume
there is a density a density $f^*(e)$ such that
%
\begin{equation}
\label{homoscedastic} f(y_1\mid x_1,x_2) = f^*
\bigl(y_1 - m\bigl(x^{\bar{m}}_1,x^{\bar{m}}_2
\bigr)\bigr)
\end{equation}
that is, the residuals all have the same distribution.
In this case, we can set $(x^{\bar{m}}_1,x^{\bar{m}}_2)$ to be all
the variables and remove $(x^{\bar{g}}_1,x^{\bar{g}}_2)$. If we set
both $x^{\bar{g}}$ and $x^{\bar{m}}$ to be the entire set $x$ we
arrive at a centered conditional density estimate
\[
f^c(\varepsilon,y_2\mid x_1,x_2) =
f\bigl(\varepsilon+m(x_1,x_2),y_2\mid
x_1,x_2\bigr).
\]
This centering can improve the finite sample performance of our
estimator at or near the homoscedastic case in which $f^c$ is close to
constant in $x_1$ and hence incurs lower bias than the uncentered version.

Here we will formalize the partition of the covariate space into
components associated with centering $y_1$ and with conditioning. To do
this, we divide $x = (x_1,x_2)$ into $(x^m,x^s,x^g)$ where $x^s$ are
the components common to both $x^{\bar{m}}= (x^m,x^s)$ and $x^{\bar
{g}}= (x^s,x^g)$ with $x^m$ and $x^g$ containing those components only
appearing one or other of the centering and conditioning variables. We
define these variables to take values on spaces $\calX^a = \real
^{d_{xa}} \otimes\, S_{xa}$ for $a \in(m,s,g)$ with $\calX= \mathcal
{X}^m\otimes\mathcal{X}^s\otimes\mathcal{X}^g$ and $\mathcal
{X}^{\bar{m}}= \mathcal{X}^m\otimes\mathcal{X}^s$ and $\mathcal
{X}^{\bar{g}}=\mathcal{X}^s\otimes\mathcal{X}^g$, similarly the
distribution of observations on these spaces will be given by
$h^a(x^a_1,x^a_2)$ for $a$ replaced by any of $(m,s,g,\bar{m},\bar{g})$.

We note that when $y_1$ is vector valued, it is not necessary to center
all of its components. The results below also encompass the case where
only some components are centered by interpreting $m(x^{\bar
{m}}_1,x^{\bar{m}}_2) = 0$ for the non-centered components. It is also
possible to include $y_2$ within $x^m_2$ (but not within $x^{\bar
{g}}_2$) without affecting these results.

The following regularity structures may be assumed in the theorems below:
\begin{enumerate}[(D4)]
\item[(D1)] $g$ is bounded and continuous in $x_1$ and $y_1$.

\item[(D2)] $\int y_1^2 g(x_1,x_2,y_1,y_2) \,\mathrm{d}y_1 < \infty$ for
all $x \in
\calX$.

\item[(D3)] All third derivatives of $g$ with respect to $x_1$ and $y_1$
exist, are
continuous and bounded.

\item[(D4)] The support of $x$, $\calX$ is compact and $h(x_1,x_2)$ is
bounded away from zero with infimum
\[
h^- = \inf_{(x_1,x_2) \in\mathcal{X}} h(x_1,x_2) > 0.
\]

\item[(D5)] The expected value function $m(x_1,x_2)$ is bounded, as is its
gradient $\nabla_{x_1} m(x_1,x_2)$.
\end{enumerate}
We note that under these conditions, continuity of $h$ and $f$ in $x_1$ and
$y_1$ is inherited from $g$. We also have that $\calX^a$ is compact
for $a \in(m,s,g,\bar{m},\bar{g})$ and similarly $h^a(x^a_1,x^a_2) >
h^-$. Assumption \textup{(D4)} is generally employed for models
involving non-parametric smoothing and is required for the uniform
convergence results that we establish; in practice it is often possible
to bound the range of values that a covariate can take. This assumption
is, however, more restrictive than required for general regression
problems and can, in fact, be removed in special cases of the methods
studied here. We have noted where this is possible below, with results
provided in Supplemental Appendix~E (Hooker \cite{Hooker14}).

In the case of centered densities (i.e., $x^{\bar{m}}$ is not
trivial), we also assume that $f$ is differentiable in $y_1$ and has a
finite second moment, uniformly over $x$:
\begin{enumerate}[(E1)]
\item[(E1)] $\sup_{(x_1,x_2) \in\calX} \sum_{y_2 \in S_y} \int
| \nabla
_{y_1} f(y_1,y_2\mid x_1,x_2) | \,\mathrm{d}y_1 <
\infty$,

\item[(E2)] $\sup_{(x^m_1,x^m_2) \in\mathcal{X}^m} \sum_{y_2 \in S_y}
\int |y_1^2 f(y_1,y_2\mid x_1,x_2) \mid \,\mathrm{d}y_1 < \infty$
\end{enumerate}
and note that these conditions need only apply to those components of
$y_1$ which are centered.

\subsection{Kernel estimators} \label{seckernels}

In order to apply the disparity methods described above, we will need
estimates of  $f^c(\varepsilon,y_2\mid x^{\bar{g}}_1,x^{\bar{g}}_2)$ which
we will obtain through kernel density and Nadaraya--Watson estimators.
Specifically, we first estimate the density of the centering variables
$(X^{\bar{m}}_1,X^{\bar{m}}_2)$:
%
\begin{equation}
\hat{h}_n^{m}\bigl(x^{\bar{m}}_1,x^{\bar{m}}_2,
\omega\bigr) = \frac
{1}{nc_{nx^{\bar{m}}_2}^{d_{x\bar{m}}}} \sum_{i=1}^n
K_x^m \biggl( \frac{x^{\bar{m}}_1 - X^{\bar{m}}_{i1}(\omega
)}{c_{nx^{\bar
{m}}_2}} \biggr) I_{x^{\bar{m}}_2}
\bigl(X^{\bar{m}}_{i2}(\omega)\bigr) \label{hathm}
\end{equation}
and define a Nadaraya--Watson estimator for the continuous response
variables $y_1$ based on them:
%
\begin{equation}
\hatmnh\bigl(x^{\bar{m}}_1,x^{\bar{m}}_2,
\omega\bigr) = \frac{ ({1}/{n
c_{nx^{\bar{m}}_2}^{d_{x\bar{m}}}})\sum_{i=1}^n Y_{i1}(\omega)
K_x^m (({x^{\bar{m}}_1-X^{\bar{m}}_{i1}(\omega
)})/{c_{nx^{\bar{m}}_2}} )
I_{x^{\bar{m}}_2} (X^{\bar{m}}_{i2}(\omega) ) }{ \hat
{h}_n^{m}(x^{\bar{m}}_1,x^{\bar{m}}_2) }. \label{hatms}
\end{equation}
We then obtain residuals from this estimator
%
\begin{equation}
\tilde{E}_i(\tilde{m},\omega) = Y_i(\omega) -
\tilde{m}\bigl(X^{\bar
{m}}_{i1}(\omega),X^{\bar{m}}_{i2}(
\omega)\bigr),\qquad i = 1,\ldots,n \label{errt}
\end{equation}
and use these with the $Y_{i2}$ to obtain a joint density estimate with
the $(X^{\bar{g}}_{i1},X^{\bar{g}}_{i2})$:
%
\begin{eqnarray}\label{hatgs}
&& \hatgnh\bigl(x^{\bar{g}}_1,x^{\bar{g}}_2,e,y_2,
\tilde{m},\omega\bigr)
\nonumber\\[-8pt]\\[-8pt]\nonumber
&&\quad  = \frac{1}{nc_{nx^{\bar{g}}_2}^{d_{x\bar{g}}}
c_{ny_2}^{d_y}} \sum_{i=1}^n
K_x \biggl( \frac{x^{\bar{g}}- X^{\bar{g}}_{i1}(\omega
)}{c_{nx^{\bar
{g}}_2}} \biggr) K_y \biggl(
\frac{e -
\tilde{E}_i(\tilde{m},\omega)}{c_{ny_2}} \biggr) I_{x^{\bar
{g}}_2}\bigl(X^{\bar{g}}_{i2}(
\omega)\bigr) I_{y_2}\bigl(Y_{i2}(\omega)\bigr).
\nonumber
\end{eqnarray}
We then estimate the density of the $(X^{\bar{g}}_{i1},X^{\bar
{g}}_{i2})$ alone
%
\begin{equation}
\hathnh\bigl(x^{\bar{g}}_1,x^{\bar{g}}_2,
\omega\bigr) = \frac
{1}{nc_{nx^{\bar{g}}_2}^{d_{x\bar{g}}}} \sum_{i=1}^n
K_x \biggl( \frac{x^{\bar{g}}_1 -
X^{\bar{g}}_{i1}(\omega)}{c_{nx^{\bar{g}}_2}} \biggr) I_{x^{\bar
{g}}_2}
\bigl(X^{\bar{g}}_{i2}(\omega)\bigr) \label{haths}
\end{equation}
and use these to obtain an estimate of the conditional distribution of
the centered responses:
%
\begin{equation}
\hatfnh(e,y_2\mid x_1,x_2,\omega) =
\frac{\hatgnh(x^{\bar{g}}_1,x^{\bar{g}}_2,e,y_2,\hatmnh,\omega
)}{\hathnh(x^{\bar{g}}_1,x^{\bar{g}}_2,\omega)} \label{hatfs}.
\end{equation}
Finally, we shift $\hatfnh$ by $\hatmnh$ to remove the centering:
%
\begin{equation}
\breve{f}_n(y_1,y_2\mid
x_1,x_2,\omega) = \hatfnh\bigl(y_1-\hatmnh
\bigl(x^{\bar{m}}_1,x^{\bar{m}}_2,\omega
\bigr),y_2\mid x^{\bar{g}}_1,x^{\bar
{g}}_2,
\omega\bigr). \label{brevef}
\end{equation}
Throughout the above, $I_{x}(X)$ is the indicator function of $X = x$
and $K_x$, $K_x^m$ and $K_y$ are densities on the spaces $\real
^{d_{x\bar{g}}
}$, $\real^{d_{x\bar{m}}}$ and $\real^{d_y}$, respectively. We have used
$c_{nx^{\bar{m}}_2}$, $c_{nx^{\bar{g}}_2}$ and $c_{ny_2}$ to
distinguish the different rates which these bandwidths will need to
follow. Further conditions on these are detailed below.

Here we have employed the errors $\tilde{E}_i(\tilde{m},\omega)$ for
the sake of notational compactness. We have defined centering by a
generic $\tilde{m}$ in (\ref{errt})--(\ref{hatgs}), which we will
employ in developing its $L_1$ convergence below, but have replaced
this with $\hatmnh$ in (\ref{hatfs}) and (\ref{brevef}) to indicate
real-world practice.

In the case of uncentered conditional density estimates ($x^{\bar{m}}$
trivial), these reduce to
%
\begin{eqnarray}
\hat{g}_n^{*}(x_1,x_2,y_1,y_2,
\omega) & =& \frac
{1}{nc_{nx_2}^{d_x}c_{ny_2}^{d_y}} \sum_{i=1}^n
K_x \biggl( \frac{x_1 - X_{i1}(\omega)}{c_{nx_2}} \biggr) K_y
\biggl(
\frac{y_1 -
Y_{i1}(\omega)}{c_{ny_2}} \biggr)
\nonumber\\[-8pt]\label{hatg} \\[-8pt]\nonumber
&&\hspace*{55pt}{}\times  I_{x_2}\bigl(X_{i2}(\omega)\bigr)
I_{y_2}\bigl(Y_{i2}(\omega)\bigr),
\\
\hat{h}_n^{*}(x_1,x_2,\omega) & =& \frac{1}{nc_{nx_2}^{d_x}} \sum_{i=1}^n
K_x \biggl( \frac{x -
X_{i1}(\omega)}{c_{nx_2}} \biggr) I_{x_2}
\bigl(X_{i2}(\omega)\bigr)
\nonumber\\[-8pt]\label{hath} \\[-8pt]\nonumber
& =& \sum_{y_2
\in S_y} \int_{\real^{d_y}}
\hat{g}_n^{*}(x_1,x_2,y_1,y_2,
\omega) \,\mathrm{d}y_1,
\nonumber
\\
\hat{f}_n^{*}(y_1,y_2\mid
x_1,x_2,\omega) & =& \frac{\hat{g}_n^{*}(x_1,x_2,y_1,y_2,\omega
)}{\hat
{h}_n^{*}(x_1,x_2,\omega)} \label{hatf}.
\end{eqnarray}
And for homoscedastic regression estimators ($x^{\bar{g}}$ and $y_2$
empty), we have
%
\begin{eqnarray}
\hatmnh(x_1,x_2,\omega) & =& \frac{ \sum_{i=1}^n Y_{i1}(\omega)
K_x (({x_1-X_{i1}(\omega)})/{c_{nx_2}} )
I_{x_2} (X_{i2}(\omega) ) }{ \sum_{i=1}^n
K_x (({x_1-X_{i1}(\omega)})/{c_{nx_2}} )
I_{x_2} (X_{i2}(\omega) )},
\label{hatm}
\\
\hat{f}_n^c(e,\omega) & =& \frac{1}{n c_{ny_2}^{d_y}} \sum
_{i=1}^n K_y \biggl(
\frac{e - (Y_i(\omega) - \hatmnh(X_{i1}(\omega),X_{i2}(\omega
)))}{c_{ny_2}} \biggr), \label{hatfss}
\\
\tilde{f}_n(y_1\mid x_1,x_2,
\omega) & =& \hat{f}_n^c \bigl(y_1-
\hatmnh(x_1,x_2,\omega),\omega\bigr) \label{hatft}
\end{eqnarray}
with notation $c_{ny_2}$ maintained as a bandwidth for the sake of consistency.

We note that while these estimates do require some extra computational
work, they are not, in fact, more computationally burdensome than the
methods proposed for independent, univariate data in Beran \cite
{Beran77}. The evaluation cost of each of the density estimates and
non-parametric smooths above is $\mathrm{O}(n)$ operations and $\breve
{f}_n(y_1,y_2\mid x_1,x_2)$ can be evaluated in a few lines of code in the
{\texttt R} programming language. In simulations reported in
Section~\ref{secsimulations} the computing time required of our
methods exceeds
that of maximum likelihood methods by a factor of 10, and alternative
robust methods by a factor of 5, rendering them very feasible in
practical situations.

Throughout we make the following assumptions on the kernels $K_x$, $K_x^m$,
and $K_y$. These will all conform to conditions on a general kernel
$K(z)$ over a Euclidean space of appropriate dimension~$\real^{d_z}$:
\begin{enumerate}[(K3)]
\item[(K1)] $K(z)$, is a density on $\real^{d_z}$.

\item[(K2)] For some finite $K^+$, $\sup_{z \in\real^{d_z}} K(z) < K^+$.

\item[(K3)]$\lim\llVert z\rrVert ^{2d_z} K(z) \rightarrow0$ as $\llVert
z\rrVert \rightarrow
\infty$.
\item[(K4)]$K(z) = K(-z)$.

\item[(K5)]$\int\llVert z\rrVert ^2 K(z) \,\mathrm{d}z < \infty$.

\item[(K6)]$K$ has bounded variation and finite modulus of continuity.
\end{enumerate}
We also assume that following properties of the bandwidths. These will
be given
in terms of the number of observations falling at each combination
values of
the discrete variables.
\begin{eqnarray*}
n\bigl(x^a_2\bigr) &=& \sum_{i=1}^n
I_{x^a_2}\bigl(X^a_{2i}(\omega)\bigr),
\qquad n(y_2) = \sum_{i=1}^n
I_{y_2}\bigl(Y_{2i}(\omega)\bigr),
\\
n \bigl(x^a_2,y_2\bigr) &=& \sum
_{i=1}^n I_{x^a_2}\bigl(X^a_{2i}(
\omega)\bigr)I_{y_2}\bigl(Y_{2i}(\omega)\bigr),
\end{eqnarray*}
where these rates are defined for $a$ covering any of $(m,s,g,\bar
{m},\bar{g})$ or the whole space. As \mbox{$n \rightarrow\infty$}:
\begin{enumerate}[(B1)]
\item[(B1)] $c_{nx_2} \rightarrow0$, $c_{ny_2} \rightarrow0$.

\item[(B2)] $n(x^a_2) c_{nx^a_2}^{d_{xa}} \rightarrow\infty$ for all $x_2
\in S_x$ and $n(x^a_2,y_2) c_{nx^a_2}^{d_{xa}} c_{ny_2}^{d_y} \rightarrow
\infty$ for all
$(x^a_2,y_2) \in S_{x^a} \otimes S_y$. \vspace*{1pt}

\item[(B3)]$ n(x^a_2) c_{nx^a_2}^{2 d_{xa}} \rightarrow\infty$.

\item[(B4)]$n(x^a_2,y_2) c_{nx^a_2}^{2d_{xa}} c_{ny_2}^{2d_y}
\rightarrow\infty$.

\item[(B5)]$\sum_{n(x^a_2)=1}^{\infty} c_{nx^a_2}^{-d_{xa}}\mathrm{e}^{-\gamma
n(x^a_2) c_{nx^a_2}^{d_x}} \leq
\infty$ for all $\gamma> 0$.

\item[(B6)]$n(y_2) c_{ny_2}^4 \rightarrow0$ if $d_y = 1$ and
$n(x^a_2) c_{nx^a_2}^4 \rightarrow0$ if $d_{xa} = 1$,
\end{enumerate}
where the sum is taken to be over all observations in the case that
$X^a_2$ or
$Y_2$ are singletons.

\subsection{Notational conventions} \label{secnotation}

Because of the complexity involved in dealing with two partitions, $x =
(x^m,x^s,x^g)$ and $x = (x_1,x_2)$, along with kernel estimators and
integrals, this paper will take some notational shortcuts; which ones
we take will differ between sections. These will allow us to ignore
notational complexities that do not affect the particular results being
discussed. Here we will forecast these.

Section~\ref{consistency1} demonstrates the consistency of
kernel-based conditional density estimates. This section will require
the distinction between continuous-valued and discrete-valued
components of $x$ and $y$ and we will emphasize the division $x =
(x_1,x_2)$. However the particular division between centering and
conditioning variables will not be important in our calculations and we
will thus suppress this notation. Formally, our results will apply to
the case where both $x^{\bar{m}}$ and $x^{\bar{g}}$ contain all the
components of $x$. However, they extend to any partition following
modification of the bandwidth scaling to reflect the dimension of the
real-valued components $(x^m,x^s,x^g)$. We have kept the notation of
$X_{i1}(\omega)$ depending on $\omega$ throughout this section
facilitate the precise description of convergence results.

In Sections~\ref{secconsistentdisparity} and~\ref{secCLT}, the
opposite case will be true. We will suppress the distinction between
discrete and continuous random variables but the partition of the
covariates into centering and conditioning components will have a
substantial effect on our results. Here, for the sake of notational
compactness we define a measure $\nu$ over $\real^{d_y} \otimes\, S_y$
and $\mu$ over $\real^{d_x} \otimes\, S_x$ given by the product of
counting and Lebesgue measure. Where needed, we will write for any
function $F(x_1,x_2,y_1,y_2)$,
%
\begin{equation}
\label{notationshortcut} \sum_{x \in S_x,y \in S_y} \int\!\!\!\int
F(x_1,x_2,y_1,y_2)
\,\mathrm{d}x_1 \,\mathrm{d}y_1 = \int\!\!\!\int F(x,y) \,\mathrm{d}\nu
(y) \,\mathrm{d} \mu(x).
\end{equation}
We\vspace*{1pt} will similarly define measures $\mu^g$, $\mu^m$, $\mu^{\bar{g}}$
and $\mu^{\bar{m}}$ over $\mathcal{X}^g$, $\mathcal{X}^m$,
$\mathcal{X}^{\bar{g}}$ and $\mathcal{X}^{\bar{m}}$, respectively.
In some places, we will refer to the centered $\varepsilon= y -
m(x^{\bar
{m}})$ where we will understand $m(x^{\bar{m}})$ to be zero on the
discrete-valued components of $y$ as well as those components of $y_1$
which are not being centered. In this context, we will subsume the
indicator functions used above within the kernel and understand
\[
K_x \biggl(\frac{x^{\bar{g}}- X^{\bar{g}}_{i}}{c_{nx^{\bar{g}}_2}}
\biggr) = K_x \biggl(
\frac{x^{\bar{g}}_1 - X^{\bar{g}}_{i1}}{\cnxgb
} \biggr)I_{x^{\bar{g}}_2}\bigl(X^{\bar{g}}_{i2}
\bigr).
\]
Here we have changed bandwidth notation to $\cnxgb$ in favor of
$c_{nx^{\bar{g}}_2}$ and understand that $c_{na}$ can depend on
$x_2^a$, but we have maintained the distinction as to which of $\bar
{m}$ or $\bar{g}$ $a$ belongs to.
We will also encounter a change of variables written as
\[
\int F\bigl(x^{\bar{g}},y\bigr)\frac{1}{\cnxgb^{d_{x\bar{g}}}}K_x
\biggl(
\frac
{x^{\bar{g}}- X^{\bar{g}}_{i}}{\cnxgb} \biggr) \,\mathrm{d}\mu
^{\bar
{g}}\bigl(x^{\bar{g}}\bigr) = \int
F\bigl(X^{\bar{g}}_i + \cnxgb u,y\bigr) K_x(u) \,\mathrm{d}u
\]
in which we will interpret $u$ as being a vector which is non-zero only
on the continuous components of $x^{\bar{g}}$. Similar conventions
will be employed for all other components of $x$ and of $y$. In these
sections, we will drop $\omega$ from our notation for the sake of
compactness and because it will be less relevant to defining our results.

\section{Consistency results for conditional densities over spaces of mixed types} \label{consistency1}

In this section, we will provide a number of $L_1$ consistency results for
kernel estimates of densities and conditional densities of multivariate random
variables in which some coordinates take values in Euclidean space
while others
take values on a discrete set. Pointwise consistency of conditional density
estimates of this form can be found in, for example,
Li and Racine \cite{LiRacine07} and Hansen \cite{Hansen04}.
However, we are unaware of equivalent $L_1$ results which will be
necessary for
our development of conditional disparity-based inference. Throughout,
we have
assumed that both the conditioning variable $x$ and the response $y$ are
multivariate with both types of coordinates. The specification to univariate
models, or models with only discrete or only continuous variables in
either $x$
or $y$ (and to unconditional densities) is readily seen to be covered
by our
results as well.

As a further generalization of the results in Li and Racine \cite{LiRacine07},
we include the centered
version of conditional density estimates defined by (\ref
{hatgs})--(\ref{brevef}). We will
demonstrate the consistency of results for these estimates, from which
consistency for
uncentered conditional densities and results for homoscedastic
conditional densities (\ref{homoscedastic}) are special cases.

Supplemental Appendix~B (Hooker \cite{Hooker14})
provides a set of intermediate results on the uniform and $L_1$
convergence of non-parametric regression and centered density estimates
of missed types. Following these, we are able to establish the uniform
(in $x$) $L_1$ (in $y$) convergence of multivariate densities:

%
\begin{teo} \label{jointsupL1}
Let $\{(X_{n1},X_{n2},Y_{n1},Y_{n2}), n \geq1\}$ be given as in
Section~\ref{secassumptions} under assumptions \textup{(D1)}--\textup{(D4)}, \textup{(K1)--(K6)} and
\textup{(B1)--(B5)} then there
exists a set $B$ with $P(B)=1$ such that for all $\omega\in B$
%
\begin{equation}
\label{jointsupL12} \sup_{(x_1,x_2) \in\mathcal{X}} \sum_{y_2 \in S_y}
\int\bigl\vert\hatgnh(x_1,x_2,y_1,y_2,
\hatmnh,\omega) - g(x_1,x_2,y_1,y_2,m)
\bigr\vert\,\mathrm{d}y_1 \rightarrow0.
\end{equation}
\end{teo}

The proof of this theorem is given in Supplemental Appendix~C.2 (Hooker \cite{Hooker14}). The results above can now
be readily extended to equivalent $L_1$ results for
conditional densities. We begin by considering centered densities and
then proceed to uncenter them.

%
\begin{teo} \label{conditionalL15}
Let $\{(X_{n1},X_{n2},Y_{n1},Y_{n2}), n \geq1\}$ be given as in
Section~\ref{secassumptions} under assumptions \textup{(D1)}--\textup{(D4)}, \textup{(K1)--(K6)} and
\textup{(B1)--(B5)}:
\begin{enumerate}
\item There exists a set $B_I$ with $P(B_I) = 1$ such that for all
$\omega\in B_I$,
%
\begin{equation}
\label{condmargoverx5} \sum_{x_2 \in S_x} \sum
_{y_2 \in S_y} \int h(x_1,x_2) \bigl\vert
\hatfnh(\varepsilon,y_2\mid x_1,x_2,\omega) -
f^c(\varepsilon,y_2\mid x_1,x_2)
\bigr\vert\,\mathrm{d}\varepsilon\,\mathrm{d}x_1 \rightarrow0.
\end{equation}
\item If further, assumptions \textup{(D4)} and \textup{(B5)}
hold, there exists a set $B_S$ with $P(B_S) = 1$ such that for all
$\omega\in B_S$:
%
\begin{equation}
\label{condsup5} \sup_{(x_1,x_2) \in\mathcal{X}} \sum_{y_2 \in S_y}
\int\bigl\vert\hatfnh(\varepsilon,y_2\mid x_1,x_2,
\omega) - f^c(\varepsilon,y_2\mid x_1,x_2)
\bigr\vert\,\mathrm{d}\varepsilon\rightarrow0.
\end{equation}
\end{enumerate}
\end{teo}

The proof of this theorem is given in Supplemental Appendix~C.2 (Hooker \cite{Hooker14}). From here, we can
examine the behavior of $\breve{f}_n$.

\begin{teo} \label{conditionalL1}
Let $\{(X_{n1},X_{n2},Y_{n1},Y_{n2}), n \geq1\}$ be given as in
Section~\ref{secassumptions} under assumptions \textup{(E1)--(E2)}, \textup{(D1)--(D4)}, \textup{(K1)--(K6)} and
\textup{(B1)--(B5)}:
\begin{enumerate}
\item There exists a set $B_I$ with $P(B_I) = 1$ such that for all
$\omega\in B_I$,
%
\begin{equation}
\label{condmargoverx} \sum_{x_2 \in S_x} \sum
_{y_2 \in S_y} \int h(x_1,x_2) \bigl\vert
\breve{f}_n(y_1,y_2\mid x_1,x_2,
\omega) - f(y_1,y_2\mid x_1,x_2)
\bigr\vert\,\mathrm{d}y_1 \,\mathrm{d}x_1 \rightarrow0.
\end{equation}
\item If further, assumptions \textup{(D4)} and \textup{(B5)}
hold, there exists a set $B_S$ with $P(B_S) = 1$ such that for all
$\omega\in B_S$:
%
\begin{equation}
\label{condsup} \sup_{(x_1,x_2) \in\mathcal{X}} \sum_{y_2 \in S_y}
\int\bigl\vert\breve{f}_n(y_1,y_2\mid
x_1,x_2,\omega) - f(y_1,y_2\mid
x_1,x_2) \bigr\vert\,\mathrm{d}y_1 \rightarrow0.
\end{equation}
\end{enumerate}
\end{teo}

\begin{pf}
We begin by writing
\begin{eqnarray*}
&& \sum_{y_2 \in S_y} \int\bigl\vert\breve{f}_n(y_1,y_2
\mid x_1,x_2,\omega) - f(y_1,y_2
\mid x_1,x_2) \bigr\vert\,\mathrm{d}y_1
\\
&&\quad \leq\sum_{y_2 \in S_y} \int\bigl\vert
\breve{f}_n(y_1,y_2\mid x_1,x_2,
\omega) - f^c\bigl(y_1 - \hatmnh(x_1,x_2),y_2
\mid x_1,x_2\bigr) \bigr\vert\,\mathrm{d}y_1
\\
&&\qquad{} + \sum_{y_2 \in S_y} \int\bigl
\vert f^c\bigl(y_1 - \hatmnh(x_1,x_2),y_2
\mid x_1,x_2\bigr) - f^c\bigl(y_1
- m(x_1,x_2),y_2\mid x_1,x_2
\bigr)\bigr\vert\,\mathrm{d}y_1
\\
&&\quad \leq\sum_{y_2 \in S_y} \int\bigl\vert
\breve{f}_n(y_1,y_2\mid x_1,x_2,
\omega) - f^c\bigl(y_1 - \hatmnh(x_1,x_2),y_2
\mid x_1,x_2\bigr) \bigr\vert\,\mathrm{d}y_1
\\
&&\qquad{} + \sup_{(x_1,x_2) \in\mathcal
{X}}\bigl|
\hatmnh(x_1,x_2) - m(x_1,x_2)\bigr|
\sum_{y_2 \in S_y} \int\bigl\vert\nabla_{y_1}
f^c(y_1,y_2\mid x_1,x_2)
\bigr\vert\,\mathrm{d}y_1.
\end{eqnarray*}
The first term of the last line converges almost surely from Theorem
\ref{conditionalL15} applied either marginalized over $(x_1,x_2)$ to
obtain (\ref{condmargoverx}) or after taking a supremum to obtain
(\ref{condsup}). The second term~follows from Theorem~B.2 in the Supplemental Appendix (Hooker \cite{Hooker14}) and
assumption~(E1).
\end{pf}

These results can now be applied to the more regular conditional
density estimates \mbox{(\ref{hatg})--(\ref{hatf})} and homoscedastic
conditional density estimates (\ref{hatm}--\ref{hatft}). For the
sake of completeness, we state these directly as corollaries without proof.

%
\begin{cor} \label{jointsupL1hat}
Let $\{(X_{n1},X_{n2},Y_{n1},Y_{n2}), n \geq1\}$ be given as in
Section~\ref{secassumptions} under assumptions \textup{(D1)--(D3)}, \textup{(K1)--(K6)} and
\textup{(B1)--(B2)} then:
\begin{enumerate}
\item For almost all $x = (x_1,x_2) \in\real^{d_x} \otimes\, S_x$ there
exists a
set $B_x$ with $P(B_x) = 1$ such that for all $\omega\in B_x$
%
\begin{equation}
\label{condatxhat} \sum_{y_2 \in S_y} \int\bigl\vert
\hat{f}_n^{*}(y_1,y_2\mid
x_1,x_2,\omega) - f(y_1,y_2\mid
x_1,x_2) \bigr\vert\,\mathrm{d}y_1 \rightarrow0.
\end{equation}
\item There exists a set $B_I$ with $P(B_I) = 1$ such that for all
$\omega\in B_I$,
%
\begin{equation}
\label{condmargoverxhat} \sum_{x_2 \in S_x} \sum
_{y_2 \in S_y} \int h(x_1,x_2) \bigl\vert
\hat{f}_n^{*}(y_1,y_2\mid
x_1,x_2,\omega) - f(y_1,y_2\mid
x_1,x_2) \bigr\vert\,\mathrm{d}y_1
\,\mathrm{d}x_1 \rightarrow0.
\end{equation}
\item If further, assumptions \textup{(D4)} and \textup{(B5)}
hold, there
exists a set $B_S$ with $P(B_S)=1$ such that for all $\omega\in B_S$
%
\begin{equation}
\label{jointsupL12hat} \sup_{(x_1,x_2) \in\mathcal{X}} \sum_{y_2
\in S_y}
\int\bigl\vert\hat{g}_n^{*}(x_1,x_2,y_1,y_2,
\omega) - g(x_1,x_2,y_1,y_2)
\bigr\vert\,\mathrm{d}y_1 \rightarrow0
\end{equation}
and
%
\begin{equation}
\label{condsuphat} \sup_{(x_1,x_2) \in\mathcal{X}} \sum_{y_2 \in S_y}
\int\bigl\vert\hat{f}_n^{*}(y_1,y_2
\mid x_1,x_2,\omega) - f(y_1,y_2
\mid x_1,x_2) \bigr\vert\,\mathrm{d}y_1
\rightarrow0.
\end{equation}
\end{enumerate}
\end{cor}

%
\begin{cor} \label{homol1}
Let $\{(X_{n1},X_{n2},Y_{n1}), n \geq1\}$ be given as in Section~\ref
{secassumptions} with the restriction (\ref{homoscedastic}), under
assumptions \textup{(D1)--(D4)}, \textup{(E1)--(E2)},
\textup{(K1)--(K6)}, \textup{(B1)--(B2)} and
\textup{(B5)} there exists a set $B$ with $P(B)=1$ such that for
all $\omega
\in B$
%
\begin{equation}
\label{margresidl1} \int\bigl\vert\hatfnh^c(e,\omega) -
f^c(e) \bigr\vert\,\mathrm{d}e \rightarrow0
\end{equation}
and
%
\begin{equation}
\label{homocondsup} \sup_{(x_1,x_2) \in\mathcal{X}} \int\bigl
\vert\tilde
{f}_n(y_1\mid x_1,x_2,\omega) -
f^c\bigl(y_1-m(x_1,x_2)\bigr)
\bigr\vert\,\mathrm{d}y_1 \rightarrow0.
\end{equation}
\end{cor}

The above theorems rely on the compactness of $\mathcal{X}$
(assumption \textup{(D4)}), this is necessary due to the estimate
$\hatmnh(x^{\bar{m}}_1,x^{\bar{m}}_2)$, and is necessary for uniform
convergence in $\mathcal{X}$. However, a weaker version can be given
for non-centered densities which does not require a compact support:

%
\begin{teo} \label{jointsupL1sup}
Let $\{(X_{n1},Y_{n1}), n \geq1\}$ be given as in Section~\ref
{secassumptions} under assumptions \textup{(D1)--(D3)}, \textup{(K1)--(K6)} and
\textup{(B1)--(B2)} then for almost all $x=(x_1,x_2)$
there exists a
set $B_{x}$ with $P(B_{x})=1$ such that for all $\omega\in
B_{x}$
%
\begin{equation}
\label{jointsupL11} \sum_{y_2 \in S_y} \int\bigl\vert
\hatgnh(x_1,x_2,y_1,y_2,\omega)
- g(x_1,x_2,y_1,y_2) \bigr
\vert\,\mathrm{d}y_1 \rightarrow0.
\end{equation}
\end{teo}

\begin{pf}
For (\ref{jointsupL11}), we observe that
\begin{eqnarray*}
\sum_{x_2 \in S_x} \sum_{y_2 \in S_y}
\int\!\!\!\int\bigl\vert\hatgnh(x_1,x_2,y_1,y_2,
\omega) - g(x_1,x_2,y_1,y_2)
\bigr\vert\,\mathrm{d}y_1 \,\mathrm{d}x_1 & =& \sum
_{x_2 \in S_x} \int T_n(x_1,x_2)
\,\mathrm{d}x_1
\\
& \rightarrow& 0
\end{eqnarray*}
almost surely with $T_n(x_1,x_2) > 0$, see \cite{DevroyeGyorfi85}, Chapter~3, Theorem~1. Thus
$T_n(x_1,x_2) \rightarrow0$ for almost all
$(x_1,x_2)$.
\end{pf}

In particular, we can rely on this theorem to remove assumption~\textup{(D4)} from the minimum disparity methods studied below in special
cases that employ $\hatgnh$ as a density estimate. Relevant further
results are given in Supplemental Appendix~E
(Hooker \cite{Hooker14}).

\section{Consistency of minimum disparity estimators for conditional models} \label{secconsistentdisparity}

In this section, we define minimum disparity estimators for the conditionally
specified models based on distributions and data defined in
Section~\ref{secassumptions}. For the purposes of notational
simplicity, we will ignore the distinction between continuous
and discrete random variables $X_1,X_2$ and $Y_1,Y_2$, but we will make
use of the division $x = (x^m,x^s,x^g)$ into those
covariates $x^m$ used to center the estimated density, those used to
condition, $x^g$, and those in both, $x^s$. We assume that a parametric
model has been proposed for these data of the form
\[
f(y\mid x) = \phi(y\mid x,\theta),
\]
where we assume that the $X_i$ are independently drawn from a
distribution $h(x)$ which is not parametrically specified. For this
model, the maximum likelihood estimator for $\theta$ given observations
$(Y_i,X_i)$, $i = 1,\ldots,n$ is
\[
\hat{\theta}_{\mathrm{MLE}} = \argmax\sum_{i=1}^n
\log\phi(Y_i\mid X_i,\theta)
\]
with attendant asymptotic variance
\[
I(\theta_0) = n \int\!\!\!\int\nabla^2_\theta\bigl[
\log\phi(y\mid x,\theta_0) \bigr] \phi(y\mid x,\theta_0)
h(x) \,\mathrm{d}\nu(y) \,\mathrm{d}\mu(x)
\]
when the specified parametric model is correct at $\theta= \theta_0$.

In the context of disparity estimation, for every value $x$ we define
the conditional disparity between $f$ and $\phi$ as
\[
D(f,\phi\mid x,\theta) = \int C \biggl( \frac{f(y\mid x)}{\phi
(y\mid x,\theta)} -1 \biggr)\phi(y\mid
x,\theta) \,\mathrm{d}\nu(y)
\]
in which $C$ is a
strictly convex function from $\real$ to $[-1\ \infty)$ with a
unique minimum
at $0$. Classical choices of $C$ include $\mathrm{e}^{-x}-1$, resulting in the
negative exponential disparity (NED) and $ [ \sqrt{x+1}-1 ]^2-1$, which
corresponds to Hellinger distance (HD).

These disparities are combined over observed $X_i$ by averaging the
disparity between $f$ and $\phi$ evaluated at each $X_i$
\[
D_n(f,\theta) = \frac{1}{n} \sum_{i=1}^n
D(f,\phi\mid X_i,\theta)
\]
(note that the $Y_i$ only appear here when $f$ is replaced by an
estimate $\breve{f}_n$)
or by integrating over the estimated density of $x^{\bar{g}}$:
\[
\tilde{D}_n(f,\theta) = \frac{1}{n} \sum
_{i=1}^n \int D\bigl(f,\phi\mid
X^m_i,x^{\bar{g}},\theta\bigr) \hathnh
\bigl(x^{\bar{g}}\bigr) \,\mathrm{d}\mu^{\bar{g}}\bigl(x^{\bar
{g}}\bigr)
\]
with limiting cases
\[
D_{\infty}(f,\theta) = \int D(f,\phi\mid x,\theta) h(x_1,x_2)
\,\mathrm{d}\mu(x)
\]
and
\[
\tilde{D}_{\infty}(f,\theta) = \int\!\!\!\int D\bigl(f,\phi\mid
x^m,x^{\bar
{g}},\theta\bigr) h^m
\bigl(x^m\bigr) h^{\bar{g}}\bigl(x^{\bar{g}}\bigr) \,\mathrm{d}
\mu^m\bigl(x^m\bigr) \,\mathrm{d}\mu^{\bar
{g}}
\bigl(x^{\bar{g}}\bigr).
\]

We now define the corresponding conditional minimum disparity estimators:
\[
\hat{\theta}^D_n = \argmin_{\theta\in\Theta}
D_n(\breve{f}_n,\theta),\qquad \tilde{\theta}^D_n
= \argmin_{\theta\in\Theta} \tilde{D}_n(\breve{f}_n,
\theta).
\]
Here we note that when the model is correct -- that is $f(y\mid x) =
\phi
(y\mid x,\theta_0)$ -- we have
that $\theta_0$ minimizes both $D_{\infty}(f,\theta)$ and $\tilde
{D}_{\infty}(f,\theta)$.

Under this definition, we first establish the existence and consistency of
$\hat{\theta}^D_n$. To do so, we note that disparity results all rely
on the
boundedness of $D(f,\phi\mid X_i,\theta)$ over $\theta$ and $f$ and a
condition of
the form that for any conditional densities $f_1$ and $f_2$,
%
\begin{equation}
\label{l1bound} \sup_{\theta\in\Theta} \bigl\vert D(f_1,
\phi\mid x,\theta) - D_n(f_2,\phi\mid x,\theta) \bigr
\vert\leq K \int\bigl\vert f_1(y\mid x) - f_2(y\mid
x) \bigr\vert\,\mathrm{d}\nu(y)
\end{equation}
for some $K > 0$. In the case of Hellinger distance (Beran \cite{Beran77}), $D(g,\theta) < 2$ and
(\ref{l1bound}) follows from Minkowski's inequality. For the alternate
class of
divergences studied in Park and Basu \cite{ParkBasu04}, boundedness of $D$ is
established
from assuming that $\sup_{t \in[-1, \infty)} |C'(t)| \leq C^* <
\infty$ which also provides
\begin{eqnarray*}
&& \biggl\vert\int\biggl[ C \biggl( \frac{f_1(y\mid x)}{\phi
(y\mid x,\theta)} - 1 \biggr) - C
\biggl( \frac{f_2(y\mid x)}{\phi(y\mid x,\theta)} - 1 \biggr)
\biggr] \phi(y\mid x,\theta) \,\mathrm{d}\nu(y) \biggr
\vert
\\
&&\quad \leq C^* \int\biggl\vert\frac{f_1(y\mid x)}{\phi(y\mid x,\theta
)} -
\frac{f_2(y\mid x)}{\phi(y\mid x,\theta)} \biggr\vert\phi(y\mid
x,\theta) \,\mathrm{d}\nu(y)
\\
&&\quad  = C^* \int\bigl\vert f_1(y\mid x) -
f_2(y\mid x) \bigr\vert\,\mathrm{d}\nu(y).
\end{eqnarray*}
For simplicity, we therefore use (\ref{l1bound}) as a condition below.

In general, we will require the following assumptions:
\begin{enumerate}[(P1)]

\item[(P1)] There exists $N$ such that  $\max_{i \in 1,\ldots,n} |\sum_{i=1}^n \phi(y|X_i,\theta_1) - \phi(y_i|X_i,\theta_2)| > 0$
with probability 1 on a nonzero set of dominating measure in $y$ whenever $n > N$ and $\theta_1 \neq \theta_2$.

\item[(P2)]$\phi(y\mid x,\theta)$ is continuous in $\theta$ for almost
every $(x,y)$.

\item[(P3)]$D_n(f,\phi\mid x,\theta)$ is uniformly bounded over $f$ in the
space of conditional densities, $(x_1,x_2) \in\mathcal{X}$ and
$\theta\in\Theta$ and (\ref{l1bound}) holds.

\item[(P4)] For every $f$, there exists a compact set $S_f \subset\Theta$
and $N$ such
that for $n \geq N$,
\[
\inf_{\theta\in S_f^c} D_{n}(f,\theta) > \inf
_{\theta\in S_f} D_{n}(f,\theta).
\]
\end{enumerate}

These assumptions combine those of Park and Basu \cite{ParkBasu04} for a
general class of
disparities with the identifiability condition \textup{(P4)}
which appears in \cite{Simpson87}, equation~(3.3), which relaxes the
assumption of compactness of $\Theta$; see also Cheng and Vidyashankar \cite
{ChenVidyashankar06}. Together, these provide the following results.

%
\begin{teo} \label{Dproperties}
Under assumptions \textup{(P1)}--\textup{(P4)}, define
%
\begin{equation}
\label{Tndef} T_n(f) = \argmin_{\theta\in\Theta} D_n(f,
\theta),
\end{equation}
for $n = 1,\ldots,\infty$ inclusive, then:
\begin{enumerate}[(iii)]
\item[(i)] For any $f \in\mathcal{F}$ there exists $\theta\in\Theta$
such that $T_n(f) = \theta$.

\item[(ii)] For $n \geq N$, for any $\theta$, $\theta= T_n(\phi(\cdot
\mid\cdot,\theta))$ is unique.

\item[(iii)] If $T_n(f)$ is unique and $f_m \rightarrow f$ in $L_1$ for each
$x$, then $T_n(f_m) \rightarrow T_n(f)$.
\end{enumerate}
The same results hold for
\[
\tilde{T}_n(f) = \argmin_{\theta\in\Theta} \tilde{D}_n(f,
\theta).
\]
\end{teo}

\begin{pf}
\textup{(i)} Existence. We first observe that it is sufficient to
restrict the infimum in (\ref{Tndef}) to $S_f$. Let $\{\theta_m\dvt
\theta_m \in
S_f\}$ be a sequence such that $\theta_m \rightarrow\theta$ as $m
\rightarrow\infty$. Since
\[
C \biggl(\frac{f(y\mid x)}{\phi(y\mid x,\theta_m)}-1 \biggr)\phi
(y\mid x,\theta_m)
\rightarrow C \biggl(\frac{f(y\mid x)}{\phi(y\mid x,\theta)}-1
\biggr)\phi(y\mid x,\theta)
\]
by assumption \textup{(P2)}, using the bound on $D(f,\phi,\theta
)$ from
assumption \textup{(P3)} we have $D_n(f,\theta_m) \rightarrow
D_n(f,\theta)$
by the dominated convergence theorem. Hence $D_n(f,t)$ is continuous in
$t$ and
achieves its minimum for $t \in S_f$ since $S_f$ is compact.

\textup{\phantom{i}(ii)} Uniqueness. This is a consequence of assumption
\textup{(P1)} and the unique minimum of $C$ at~$0$.

\textup{(iii)} Continuity in $f$. For any sequence $f_m(\cdot\mid x)
\rightarrow f(\cdot\mid x)$ in $L_1$ for every $x$ as $m \rightarrow
\infty$,
we have
%
\begin{equation}
\label{Dequicontinuity} \sup_{\theta\in\Theta} \bigl\vert D_n(f_m,
\theta) - D_n(f,\theta) \bigr\vert\rightarrow0
\end{equation}
from assumption \textup{(P3)}.

Now consider $\theta_m = T_n(f_m)$. We first observe that there exists
$M$ such
that for $m \geq M$, $\theta_m \in S_f$ otherwise from (\ref{Dequicontinuity})
and assumption \textup{(P4)}
\[
D_n(f_m,\theta_m) > \inf
_{\theta\in S_f} D_n(f_m,\theta)
\]
contradicting the definition of $\theta_m$.

Now suppose that $\theta_m$ does not converge to $\theta_0$. By the
compactness
of $S_f$ we can find a subsequence $\theta_{m'} \rightarrow\theta^*
\neq
\theta_0$ implying $D_n(f,\theta_{m'}) \rightarrow D_n(f,\theta^*)$ from
assumption \textup{(P2)}. Combining this with (\ref{Dequicontinuity})
implies $D_n(f,\theta^*) = D_n(f,\theta_0)$, contradicting the
assumption of
the uniqueness of~$T_n(f)$.
\end{pf}

%
\begin{teo} \label{consistency}
Let $\{(X_{n1},X_{n2},Y_{n1},Y_{n2}), n \geq1\}$ be given as in
Section~\ref{secassumptions} and define
\[
\theta^0_n = \argmin_{\theta\in\Theta}
D_{n}(f,\theta)
\]
for every $n$ including $\infty$. Further, assume that $\theta
^0_{\infty}$ is
unique in the sense that for every $\varepsilon$ there exists $\delta$
such that
\[
\bigl\llVert\theta- \theta^0_{\infty}\bigr\rrVert> \varepsilon\quad
\Rightarrow\quad D_{\infty
}(f,\theta) > D_{\infty}\bigl(f,
\theta^0_{\infty}\bigr) + \delta
\]
then under assumptions \textup{(D1)--(D4)}, \textup{(K1)--(K6)}, \textup{(B1)--(B2)} and
\textup{(P1)}--\textup{(P4)}:
\[
\hat{\theta}_n = T_n(\breve{f}_n)
\rightarrow\theta^0_{\infty}\qquad\mbox{as } n \rightarrow\infty\mbox{ almost surely.}
\]
Similarly,
\[
\tilde{T}_n(\breve{f}_n) = \argmin_{\theta\in\Theta}
\tilde{D}_{n}(\breve{f}_n,\theta) \rightarrow\tilde{
\theta}^0_\infty\qquad\mbox{as } n \rightarrow\infty\mbox{ almost surely.}
\]
\end{teo}

\begin{pf}
First, we observe that for every $f$, it is sufficient to restrict
attention to
$S_f$ and that
%
\begin{equation}
\label{Dnconv} \sup_{\theta\in S_f} \bigl\vert D_n(f,
\theta) - D_{\infty}(f,\theta)\bigr\vert\rightarrow0\qquad\mbox{almost surely}
\end{equation}
from the strong law of large numbers, the compactness of $S_f$ and the assumed
continuity of $C$ and of $\phi$ with respect to $\theta$.

Further,
%
\begin{eqnarray}\label{Dfnconv}
\sup_{m \in\nat,\theta\in\Theta} \bigl\vert D_m(\breve{f}_n,
\theta) - D_{m}(f,\theta)\bigr\vert& \leq& C^* \sup
_{x \in\mathcal{X}} \int\bigl\vert\breve{f}_n(y\mid x) - f(y
\mid x) \bigr\vert\,\mathrm{d}\nu(y)
\nonumber\\[-8pt]\\[-8pt]\nonumber
& \rightarrow& 0\qquad\mbox{almost surely},
\end{eqnarray}
where the convergence is obtained from Theorem~\ref{conditionalL1}.\vspace*{1pt}

Suppose that $\hat{\theta}_n$ does not converge to $\theta^0_{\infty
}$, then we
can find $\varepsilon> 0$ and a subsequence $\hat{\theta}_{n'}$ such that
$\llVert \hat{\theta}_{n'} - \theta^0_{\infty} \rrVert >
\varepsilon$ for all
$n'$. However,
on this subsequence
\begin{eqnarray*}
D_{n'}(\breve{f}_{n'},\hat{\theta}_{n'}) & =&
D_{n'}(\breve{f}_{n'},\theta_0) +
\bigl(D_{n'}(f,\theta_0)-D_{n'}(
\breve{f}_{n'},\theta_0)\bigr) + \bigl(D_\infty
(f,\theta_0) - D_{n'}(f,\theta_0)\bigr)
\\
&&{} + \bigl(D_\infty(f,\hat{\theta}_{n'}) -
D_\infty(f,\theta_0)\bigr)
\\
&&{} + \bigl(D_{n'}(f,\hat{\theta}_{n'}) -
D_\infty(f,\hat{\theta}_{n'})\bigr) + \bigl(D_{n'}(
\breve{f}_{n'},\hat{\theta}_{n'}) - D_{n'}(f,
\hat{\theta}_{n'})\bigr)
\\
& \leq& D_{n'}(\breve{f}_{n'},\theta_0) +
\delta
\\
&&{} - 2 \sup_{\theta\in\Theta} \bigl\vert
D_{n'}(f,\theta) - D(f,\theta) \bigr\vert- 2 \sup
_{\theta\in\Theta} \bigl\vert D_{n'}(\breve{f}_{n'},
\theta) - D_{n'}(f,\theta) \bigr\vert
\end{eqnarray*}
but from (\ref{Dnconv}) and (\ref{Dfnconv}) we can find $N$ so that
for ${n'}
\geq N$
\[
\sup_{\theta\in\Theta} \bigl\vert D_{n'}(f,\theta) - D(f,
\theta) \bigr\vert\leq\frac{\delta}{6}
\]
and
\[
\sup_{\theta\in\Theta} \bigl\vert D_{n'}(
\breve{f}_{n'},\theta) - D_{n'}(f,\theta) \bigr\vert\leq
\frac{\delta}{6}
\]
contradicting the optimality of $\hat{\theta}_{n'}$. The proof for
$\tilde{T}_n(\breve{f}_n)$ follows analogously.
\end{pf}

The compactness assumption \textup{(D4)} used above can be removed
for the special case of an uncentered density employed with our second
estimator: $\tilde{T}(\hat{f}_n^{*})$. This is stated in Theorem~E.1 in Supplemental Appendix~E
(Hooker \cite{Hooker14}).

\section{Asymptotic normality and efficiency of minimum disparity estimators for conditional models} \label{secCLT}

In this section, we demonstrate the asymptotic normality and efficiency of
minimum conditional disparity estimators. In order to simplify some of
our expressions, we
introduce the following notation, that for a column vector $A$ we
define the matrix
\[
A^{TT} = AA^T.
\]
This will be particularly useful in defining information matrices.

We will also frequently use the notation $y = (y_1,y_2)$ and $x =
(x_1,x_2)$, ignoring the
distinction between real and discrete valued variables. It will be
particularly relevant to distinguish $x^{\bar{g}}$ and $x^{\bar{m}}$
along with their subsets $x^g$ and $x^m$ that are solely in $x^{\bar
{g}}$ or $x^{\bar{m}}$, respectively, along with the shared dimensions
$x^s$. Because our notation would otherwise become unwieldy, we will
subsume indicator functions within kernels, and, for example, understand
\[
K_x \biggl( \frac{x^{\bar{g}}- X^{\bar{g}}_i}{\cnxgb} \biggr) =
K_x \biggl(
\frac{x^{\bar{g}}_1 - X^{\bar{g}}_{i1}}{c_{nx^{\bar
{g}}_2}} \biggr)I_{x^{\bar{g}}_2}\bigl(X^{\bar{g}}_{i2}
\bigr),
\]
where we have also suppressed the $x_2$ indicator in the bandwidth
$\cnxgb$.
Within this context, we will also occasionally abuse notation when
changing variables and write $x^{\bar{g}}+ \cnxgb v$ in which we
understand that the additive term only corresponds to the
continuous-valued entries in $x^{\bar{g}}$. We will also express
integration with respect to the distribution $\mu(x)$ and $\nu(y)$
and denote $\mu^g$, $\mu^m$, $\mu^s$, $\mu^{\bar{g}}$ and $\mu
^{\bar{m}}$ the measures marginalized to the corresponding dimensions
of $\mathcal{X}$.

The proof techniques employed
here are an extension of those developed in i.i.d. settings in
Beran \cite{Beran77}; Tamura and Boos \cite{TamuraBoos86};
Lindsay \cite{Lindsay94};
Park and Basu \cite{ParkBasu04}. In particular
we will
require the following assumptions:
\begin{enumerate}[(N3)]
\item[(N1)] Define
\[
\Psi_\theta(x,y) = \frac{\nabla_\theta
\phi(y\mid x,\theta)}{\phi(y\mid x,\theta)}
\]
then
\[
\sup_{x \in\mathcal{X}} \int\Psi_\theta(x,y)
\Psi_\theta(x,y)^T f(y\mid x) \,\mathrm{d}\nu(y) < \infty
\]
elementwise. Further, there exists $a_y > 0 $ such that
\[
\sup_{x \in\mathcal{X}} \sup_{\llVert t\rrVert \leq a_y} \sup
_{\llVert s\rrVert \leq a_y} \int\Psi_{\theta}(x_1+s,x_2,y_1+t,y_2)^2
f(y\mid x) \,\mathrm{d}\nu(y) < \infty
\]
and
\[
\sup_{x \in\mathcal{X}} \sup_{\llVert t\rrVert \leq a_y} \sup
_{\llVert s\rrVert \leq
a_y} \int\bigl(\nabla_{y_1} \Psi_\theta(x_1+t,x_2,y_1+s,y_2)
\bigr)^2 f(y\mid x) \,\mathrm{d}\nu(y)< \infty,
\]
and
\[
\sup_{x \in\mathcal{X}} \sup_{\llVert t\rrVert \leq a_y} \sup
_{\llVert s\rrVert \leq
a_y} \int\bigl(\nabla_x \Psi_\theta(x_1+t,x_2,y_1+s,y_2)
\bigr)^2 f(y\mid x) \,\mathrm{d}\nu(y) < \infty.
\]

\item[(N2)] There exists sequences $b_n$ and $\alpha_n$ diverging to
infinity along with a constant $c>0$ such that:
\begin{longlist}[(iii)]
\item[(i)] $n K_x(b_n/c_{nx}) \rightarrow0$,
$nK_y(b_n/c_{ny}) \rightarrow0$ and
\[
n \sup_{x \in\mathcal{X}} \sup_{\llVert u\rrVert > b_n} \int\!\!\!\int
_{\llVert v\rrVert
> b_n} \Psi_\theta^2(x+c_{nx}
u,y+c_{ny}v)K_y^2(u) K_x^2(v)
g(x,y) \,\mathrm{d}v \,\mathrm{d}\nu(y) \rightarrow0
\]
elementwise.

\item[(ii)] $\sup_{x \in\mathcal{X}} n P(\llVert Y_1-c_{ny} b_n\rrVert >
\alpha_n -
c) \rightarrow0$.

\item[(iii)]
\[
\sup_{x \in\mathcal{X}} \frac{1}{\sqrt{n}
c_{nx}^{d_{x}} c_{ny}^{d_{y}}} \int_{\llVert y_1\rrVert \leq\alpha
_n + c}
\bigl\vert\Psi_\theta(x,y)\bigr\vert\,\mathrm{d}\nu(y)
\rightarrow0.
\]
\item[(iv)]
\[
\sup_{x \in\mathcal{X}} \sup_{\llVert t\rrVert \leq b_n} \sup
_{\llVert s\rrVert \leq b_n} \sup_{\llVert y_1\rrVert < \alpha_n}
\frac{
g(x+c_{nx}s,y+c_{ny} t)}{g(x,y)} = \mathrm{O}(1).
\]
\end{longlist}

\item[(N3)] $\sup_{y,x} \sqrt{\phi(y\mid x,\theta)} \nabla_\theta\Psi
_\theta(y,x) = S < \infty$.

\item[(N4)] $C$ is either given by Hellinger distance $C(x) = [ \sqrt
{x+1} - 1 ]^2-1$ or
\begin{eqnarray*}
A_1(r) &=& -C''(r-1) r,\qquad
A_2(r) = C(r-1) - C'(r-1) r,
\\
A_3(r) &=& C''(r-1) r^2
\end{eqnarray*}
are all bounded in absolute value as is $r^2C^{(3)}(r)$.
\end{enumerate}

Assumption \textup{(N1)} ensures that the likelihood score
function is well controlled including for small location changes of
$x_1$ and $y_1$.
Assumption~\textup{(N2)} requires $\Psi_{\theta}$ and $y_1$ to have
well-behaved tails relative to $K_y$. In particular, assumption (N2)(i) allows us to truncate the kernels at $b_n$ which will prove
mathematically convenient throughout the remainder of the section.
Assumption (N3) concerns the regularity of the parametric model
and in particular ensures that the second derivative of Hellinger
distance with respect to parameters is well behaved. Assumption (N4)
is a restatement of conditions on the residual
adjustment function in Lindsay \cite{Lindsay94} and Park and Basu  \cite
{ParkBasu04}; a wide class of disparities satisfy these conditions
including NED, we refer the reader to Lindsay \cite{Lindsay94} for a
more complete discussion. As was the case for assumption \textup{(P3)}, we treat Hellinger distance separately in assumption (N4)
as it does not conform to the general assumptions on~$C$, but the relevant bounds can be demonstrated by other means in the
proof of Theorem~\ref{generalefficiency} below.


The demonstration of a central limit theorem involves bounding the
score function for a general disparity in terms of that for
Hellinger distance and then taking Taylor expansion of this score. For
this we need two lemmas. The first is that the weighted Hellinger distance
between $\hatfnh$ and its expectation is smaller than $\sqrt{n}$.
This is used in Theorem~\ref{generalefficiency} to remove terms
involving~$\sqrt{\hatfnh}$.

%
\begin{lem} \label{hatfHD}
Let $\{(X_{n},Y_{n}), n \geq1\}$ be given as in Section~\ref
{secassumptions}, under assumptions \textup{(D1)--(D4)}, \textup{(K1)--(K6)}, \textup{(B1)--(B4)}, and
\textup{(N1)--(N2)(iv)} for any function $J(y,x)$
satisfying the conditions on~$\Psi$ in assumptions
\textup{(N1)--(N2)(iv)}
%
\begin{equation}
\label{hatfHDrootn} \sqrt{n} \sup_{x \in\mathcal{X}} \int\!\!\!\int
J\bigl(e +
\hatmnh\bigl(x^{\bar{m}}\bigr),x\bigr) \biggl( \sqrt{\hatfnh
(e\mid x)} -
\sqrt{\frac{E
\hatgnh(x,e,\hatmnh)\mid\hatmnh}{E \hathnh(x)}} \biggr)^2 \,
\mathrm{d}\nu(e)\rightarrow0
\end{equation}
in probability and
%
\begin{eqnarray}
\label{hatfHDrootn2}
&& \frac{1}{\sqrt{n}}\sum_{i = 1}^n
\int\!\!\!\int\biggl( \sqrt{\hatfnh\bigl(e\mid X^m_i,x^g
\bigr)} - \sqrt{\frac{E
\hatgnh(X^m_i,x^g,e,\hatmnh)\mid\hatmnh}{ E\hathnh(x^g)}} \biggr
)^2
\nonumber\\[-8pt]\\[-8pt]\nonumber
&&\hspace*{48pt}{}\times  J\bigl(e + \hatmnh\bigl(X^m_i
\bigr),X^m_i,x^g\bigr)\hathnh
\bigl(x^g\bigr) \,\mathrm{d}\nu(e)\,\mathrm{d}\mu^g\bigl
(x^g\bigr)
\rightarrow0
\nonumber
\end{eqnarray}
and
%
\begin{equation}
\label{hatfHDrootn3} \frac{1}{\sqrt{n}} \sum_{i = 1}^n
\int\!\!\!\int J\bigl(e + \hatmnh\bigl(X^m_i
\bigr),X_i\bigr) \biggl( \sqrt{\hatfnh(e\mid X_i)} -
\sqrt{\frac{E
\hatgnh(X_i,e,\hatmnh)\mid\hatmnh}{\hathnh(X_i)}} \biggr)^2 \,
\mathrm{d}\nu(e)\rightarrow0.
\end{equation}
\end{lem}

The proof of this lemma is given in Supplemental Appendix~D.2 (Hooker \cite{Hooker14}).\vspace*{1pt}

A second lemma states that integrating a function $J(y,x)$ with respect
to $\hatfnh$ yields a central limit theorem. In the below, we have
used a subscript $x_*$ to help differentiate which components are being
integrated with respect to which measure.

%
\begin{lem} \label{hatfCLT}
Let $\{(X_{n},Y_{n}), n \geq1\}$ be given as in Section~\ref
{secassumptions}, under assumptions \textup{(D1)--(D4)},
\textup{(E1)--(E2)}, \textup{(K1)--(K6)},
\textup{(B1)--(B2)} and \textup{(P1)}--\textup{(P4)}
for any for any function $J(y,x)$ satisfying the conditions on
$\Psi$ in assumptions \textup{(N1)--(N2)(iv)}
and
\[
V_J = \int\!\!\!\int U_J^{TT}(
\varepsilon,x)g^c(x,\varepsilon) \,\mathrm{d}\nu(\varepsilon) \,
\mathrm{d}\mu(x) < \infty,
\]
where
\begin{eqnarray*}
U_J(\varepsilon,x) & =& \int J\bigl(\varepsilon+ m\bigl(x^m_*,x^s
\bigr),x^m_*,x^{\bar
{g}}\bigr)h^m
\bigl(x^m_*\bigr) \,\mathrm{d} \mu^m\bigl(x^m_*\bigr)
\\[-2pt]
&&{} + \int\!\!\!\int J\bigl(e_*+m\bigl(x^m,x^s_*
\bigr),x^m,x^{\bar
{g}}_*\bigr)g^c(x_*,e_*) \,\mathrm{d}\nu(e_*)
\,\mathrm{d}\mu(x_*)
\\[-2pt]
&&{} + \varepsilon\int\!\!\!\int\nabla_y J\bigl(e_*+m
\bigl(x^{\bar
{m}}\bigr),x^{\bar{m}},x^g_*\bigr)
\\[-2pt]
&&\hspace*{34pt}{}\times \frac{g^c(x^m_*,
x^s,x^g_*,e_*)h^m(x^m)}{h^{\bar{m}}(x^{\bar{m}})} \,\mathrm{d}\nu
(e_*) \,\mathrm{d}\mu^g\bigl(x^g_*\bigr) \,\mathrm{d}
\mu^m\bigl(x^m_*\bigr)
\\[-2pt]
&&{} - \varepsilon\int\!\!\!\int\nabla_y J\bigl(e_*+m
\bigl(x^m_*,x^s\bigr),x^m_*,x^s,x^g_*
\bigr)
\\[-2pt]
&&\hspace*{34pt}{}\times \frac{g^c(x^{\bar
{m}},x^g_*,e_*)h^m(x^m_*)}{h^{\bar{m}}(x^{\bar{m}})} \,\mathrm{d}\nu
(e_*) \,\mathrm{d}\mu^g\bigl(x^g_*
\bigr) \,\mathrm{d}\mu^{m}\bigl(x^m_*\bigr)
\end{eqnarray*}
elementwise, then
%
\begin{equation}
\label{tDCLT} \sqrt{n} \Biggl[ \frac{1}{n} \sum
_{i=1}^n \int\!\!\!\int J\bigl(y,X^m_i,x^{\bar{g}}
\bigr) \hatgnh(y,x) \,\mathrm{d}\nu(y)\,\mathrm{d}\mu(x) - B_n
\Biggr]\rightarrow
N(0,V_J)
\end{equation}
in distribution where
\[
B_n = 2 \int\!\!\!\int J\bigl(e + m\bigl(x^m_*,x^{\bar{g}}
\bigr),x^s_*,x^g\bigr) E \hatgnh\bigl(x^g,e,m
\bigr) h^m\bigl(x^m_*\bigr) \,\mathrm{d}\nu(e) \,\mathrm{d}\mu(x).
\]

Similarly, if
\[
\tilde{V}_J = \int\!\!\!\int\tilde{U}_J^{TT}(
\varepsilon,x)g^c(x,\varepsilon) \,\mathrm{d}\nu(\varepsilon)\,
\mathrm{d}\mu(x) < \infty,
\]
where
\begin{eqnarray*}
\tilde{U}_J(\varepsilon,x)
&=& \int\!\!\!\int J\bigl(e_* + m\bigl(x^m_*,x^s
\bigr),x\bigr) \frac
{g(x^m_*,x^{\bar{g}},e_*)}{h^{\bar{g}}(x^{\bar{g}})} \,\mathrm
{d}\mu^m\bigl(x^m_*\bigr)
\,\mathrm{d}\nu(e)
\\[-2pt]
&&{} + \int\!\!\!\int J\bigl(\varepsilon+m\bigl(x^m_*,x^s
\bigr),x^m_*,x^{\bar
{g}}\bigr)\frac{h(x^m_*,x^{\bar{g}})}{h^{\bar{g}}(x^{\bar{g}})} \,
\mathrm{d}\mu
^m\bigl(x^m_*\bigr)
\\[-2pt]
&&{} + \varepsilon\int\!\!\!\int\nabla_y J\bigl(e_*+m
\bigl(x^{\bar
{m}}\bigr),x^{\bar{m}},x^g_*\bigr)
\frac{g(x^m_*,x^s,x^g_*,e_*) h(x^m,x^{\bar
{g}}_*)}{h^{\bar{g}}(x^s,x^g_*)h^{\bar{m}}(x^{\bar{m}})} \,\mathrm
{d}\nu(e_*) \,\mathrm{d}\mu(x_*)
\\[-2pt]
&&{} - \varepsilon\int\!\!\!\int\!\!\!\int\nabla_y J\bigl(e_*+m
\bigl(x^m_*,x^s\bigr),x^m_*,x^s,x^g_*
\bigr)
\\[-2pt]
&&\hspace*{41pt}{}\times \frac{g^c(x^{\bar
{m}},x^g_*,e_*)h(x^m_*,x^s,x^g_*)}{h^{\bar{g}}(x^s,x^g_*) h^m(x^{\bar
{m}})} \,\mathrm{d}\nu(e_*) \,\mathrm{d}\mu^{m}\bigl(x^m_*
\bigr)\,\mathrm{d}\mu^g\bigl(x^g_*\bigr)
\\[-2pt]
&&{} + \int\!\!\!\int J\bigl(e_*+m\bigl(x^m_*,x^s
\bigr),x^m_*,x^g\bigr) \frac
{g^c(x^m_*,x^{\bar{g}},e_*)}{h^{\bar{g}}(x^{\bar{g}})} \,\mathrm
{d}\nu(e_*)\,\mathrm{d}
\mu^m\bigl(x^m_*\bigr)
\end{eqnarray*}
elementwise, then
%
\begin{equation}
\label{DCLT} \sqrt{n} \Biggl[ \frac{1}{n} \sum
_{i=1}^n \int J(y,X_i)
\breve{f}_n(y\mid X_i) \,\mathrm{d}\nu(y) - \tilde{B}_n
\Biggr] \rightarrow N(0,\tilde{V}_J)
\end{equation}
in distribution with
\[
\tilde{B}_n = 2 \int\!\!\!\int J\bigl(e + m\bigl(x^{\bar{m}}\bigr),x
\bigr) \frac{ E\hatgnh(x^g,e,m)}{E \hathnh(x^{\bar
{g}}) } h(x) \,\mathrm{d}\nu(e) \,\mathrm{d}\mu(x).
\]
\end{lem}

The proof of this lemma is reserved to Supplemental Appendix~D.1  (Hooker \cite{Hooker14}).

The bias and variance terms found in this lemma are rather complex due
to their generality and it will be helpful here to note the resulting
expressions for four simplifying cases and the consequence of these.
Further, in Theorem~\ref{generalefficiency} we will investigate
%
\begin{equation}
\label{eqscorereplace} \Psi_{\theta}(y\mid x) = \frac{\nabla
_{\theta} \phi_{\theta
}(y\mid x)}{\phi_{\theta}(y\mid x) },
\end{equation}
where if $\phi_{\theta}(y\mid x)$ has the form $\phi( y - m(x^{\bar
{m}};\theta)\mid x^{\bar{g}};\theta)$ we have that
\[
\Psi_{\theta} (y\mid x) = - \frac{\partial_{\theta} m(x^{\bar
{m}},\theta)\,\partial_y\phi(y-m(x^{\bar{m}};\theta
_1)| x^{\bar
{g}};\theta)}{\phi(y-m(x^{\bar{m}};\theta)|x^{\bar{g}};\theta)} +
\frac{\partial_{\theta} \phi(y-m(x^{\bar{m}};\theta)|x^{\bar
{g}};\theta)}{\phi(y-m(x^{\bar{m}};\theta)|x^{\bar{g}};\theta)},
\]
where $\partial$ is used to represent a partial gradient and $\nabla$
the total gradient. We also have that
\begin{eqnarray*}
\nabla_y J(y\mid x) & =& - \frac{\partial_{\theta} m(x^{\bar
{m}},\theta
)D^2_y \phi(y-m(x^{\bar{m}};\theta
_1)|x^{\bar{g}};\theta)}{\phi
(y-m(x^{\bar{m}};\theta)|x^{\bar{g}};\theta)}
\\
&&{}+ \frac{D^2_{\theta
y} \phi(y-m(x^{\bar{m}};\theta)|x^{\bar{g}};\theta)}{\phi
(y-m(x^{\bar{m}};\theta)|x^{\bar{g}};\theta)}
\\
&&{} + \frac{\partial_{\theta} m(x^{\bar{m}},\theta)\,
\partial_y\phi(y-m(x^{\bar{m}};\theta_1)|x^{\bar{g}};\theta
)}{\phi(y-m(x^{\bar{m}};\theta)|x^{\bar{g}};\theta)} \frac
{\partial_y \phi(y-m(x^{\bar{m}};\theta_1)|x^{\bar{g}};\theta
)^T}{\phi(y-m(x^{\bar{m}};\theta)|x^{\bar{g}};\theta)}
\\
&&{} - \frac{\partial_{\theta} \phi(y-m(x^{\bar
{m}};\theta)|x^{\bar{g}};\theta)}{\phi(y-m(x^{\bar{m}};\theta
)|x^{\bar{g}};\theta)} \frac{\partial_y \phi(y-m(x^{\bar
{m}};\theta_1)|x^{\bar{g}};\theta)^T}{\phi(y-m(x^{\bar{m}};\theta
)|x^{\bar{g}};\theta)},
\end{eqnarray*}
where we take $\partial_y^2 \phi$ to be the Hessian with respect to
$y$ and $\partial^2_{\theta y} \phi$ to be the corresponding matrix
of cross derivatives. In each of these cases, we demonstrate that
substituting in $f(y\mid x) = \phi_{\theta}(y\mid x)$ results in variance
terms given by the Fisher information
\begin{eqnarray*}
I(\theta) = \int\!\!\!\int\frac{\nabla_{\theta} \phi_{\theta
}(e\mid x)\nabla_{\theta} \phi_{\theta}(e\mid x)^T}{\phi_{\theta
}(e\mid x)} h(x) \,\mathrm{d}\nu(e)\,\mathrm{d}\mu(x)
\end{eqnarray*}
or the equivalent based on centering by $m(x,\theta)$ above.

%
%
\textit{Non-centered}: $x^{\bar{m}}= \phi$. This corresponds to the
simplest case of a conditional density estimate. Here we have
\begin{eqnarray*}
U_J(y,x) & =& J(y,x),
\\
B_n & =& 2 \int\!\!\!\int J(y,x) E \hatgnh(x,y) \,\mathrm{d}\nu(y) \,
\mathrm{d}\mu(x).
\end{eqnarray*}
We remark here that the bias $B_n$ corresponds to the bias found in
Tamura and Boos \cite{TamuraBoos86} for multivariate observations. As observed
there, the bias
in the estimate $\hatgnh$ is $\mathrm{O}(\cnxgb^2 + c_{ny}^2)$ and that of
$\hathnh$ is $\mathrm{O}(\cnxgb^2)$, regardless of the dimension of $x^{\bar
{g}}_1$ and $y_1$.
However the variance is of order $n^{-1} \cnxgb^{d_{x\bar{g}}}c_{ny}^{d_y}$
(corresponding to assumption (B2)), meaning that for
$d_x+d_y >3$,
the asymptotic bias in the Central Limit theorem is $\sqrt{n} \cnxgb^2
c_{ny}^2 \rightarrow\infty$ and will not become zero when the
variance is
controlled. We will further need to restrict to $n
\cnxgb^{2d_{x\bar{g}}}c_{ny}^{2d_y} \rightarrow\infty$, effectively
reducing the unbiased central
limit theorem to the cases where there is only one continuous variable,
although it can be either in $y$ or $x$. As in Tamura and Boos \cite
{TamuraBoos86} we also note that this bias is often small in practice;
Section~\ref{secbootstrap} demonstrates that a bootstrap method can
remove it. We also note that in this case, the assumption of a compact
domain for the covariates $x$ can be relaxed.

In the case of (\ref{DCLT}), we have
\begin{eqnarray*}
\tilde{U}_J(y,x) & =& J(y,x) + 2 \int J(y_*,x) f(y_*\mid x) \,\mathrm{d}
\nu(y_*),
\\
\tilde{B}_n & =& 2 \int\!\!\!\int J(y,x) \frac{E \hatgnh(y,x)}{E \hathnh
(x)} h(x) \,\mathrm{d}\nu(y) \,\mathrm{d}
\mu(x),
\end{eqnarray*}
where we note the additional variance due to the summation over $X_i$
values. In this case, the assignment (\ref{eqscorereplace})
with $f(y\mid x) = \phi_{\theta}(y\mid x)$ gives us that the variance
is the
information matrix directly. For $\tilde{U}_J$, we observe that
\[
\int J(y_*,x) f(y_*\mid x) \,\mathrm{d}\nu(y_*) = \int\nabla
_{\theta} \phi
_{\theta}(y\mid x)\,\mathrm{d}\nu(y_*) = 0
\]
since $\phi_{\theta}(y\mid x)$ integrates to 1 for each $x$ and each
$\theta$, yielding the same variance term as above. The bias here is
of the same order as above.

\textit{Homoscedastic}: $x^{\bar{g}}= \phi$.
Here the density estimate
assumes that $y$ has a location-scale family with $y - m(x)$
independent of $x$. In this case,
\begin{eqnarray*}
U_J(\varepsilon,x) & =& \int J\bigl(\varepsilon+m(x),x_*\bigr)
h(x_*) \,\mathrm{d}
\mu(x_*)
\\
&&{}+ \int\!\!\!\int J\bigl(e_*+m(x),x\bigr) g^c(x_*,e_*) \,\mathrm
{d}\mu(x_*) \,\mathrm{d}
\nu(e_*)
\\
&&{} + \varepsilon\int\!\!\!\int\nabla_y J\bigl(e_*+m(x),x\bigr)
g^c(x_*,e_*) \,\mathrm{d}\mu(x_*) \,\mathrm{d}\nu(e_*)
\\
&&{} - \varepsilon\int\!\!\!\int\nabla_yJ\bigl(e_* + m(x_*),x_*
\bigr) g^c(x_*,e_*) \,\mathrm{d}\nu(e) \,\mathrm{d}\mu(x_*),
\\
B_n & =& 2\int\!\!\!\int J\bigl(e+m(x),x\bigr) E\hatgnh(x,e,m) h(x) \,
\mathrm{d}\nu(e)
\,\mathrm{d}\mu(x).
\end{eqnarray*}
Here we observe that the bias is again of order $c_{nx}^2$. However,
for $e$ and $x^m$ both univariate it is possible to make $\sqrt{n} B_n
\rightarrow0$ while retaining consistency of $\hatgnh(e,m)$ and
$\hatmnh(x^m)$.

We also have
\[
\tilde{U}_J(\varepsilon,x) = U_J(\varepsilon,x),\qquad
\tilde{B}_n = B_n
\]
since in this case, both estimators are equal.

When we make the replacement (\ref{eqscorereplace}), we assume that
the assumed residual density $\phi(e;\theta)$ is parameterized so that
\[
\phi(e;\theta) = \phi^*(S_{\theta} e;\theta)
\]
with
\[
\int e\mathrm{e}^T \phi^*(e,\theta) \,\mathrm{d}\nu(e) = \int\frac{\nabla
_e \phi
^*(e,\theta)^{TT}}{\phi^*(e;\theta)} \,\mathrm{d}e = I\quad\mbox{and}\quad \int e \phi^*(e;\theta) = 0
\]
for all $\theta$ where $I$ is the $d_y\times d_y$ identity matrix. The
second equality can always be achieved by re-parameterizing so that
$\phi^*(e;\theta) = \phi( I(\theta)^{1/2} e;\theta)$
along with appropriate centering. The first equality requires that the
variance in $\phi^*(e;\theta)$ be equal to the Fisher information for
the location family $\phi^*(e+\mu;\theta)$; this condition is
satisfied, for example, for the multivariate normal density. We now
have that the total gradient is
\[
\nabla_e \phi^*(S_\theta e;\theta) = S_{\theta}
\partial_e \phi^*(S_{\theta}e;\theta)
\]
and hence
\[
U_J(\varepsilon,x) = \overline{\partial_{\theta} m}
S_{\theta} \frac
{\partial_{y} \phi^*(S_{\theta}\varepsilon;\theta)}{ \phi
^*(S_{\theta}\varepsilon;\theta)} + \frac{ \partial_{\theta} \phi
(\varepsilon;\theta)}{\phi(\varepsilon;\theta)} + \varepsilon
\bigl(
\partial_{\theta}m(x,\theta) - \overline{\partial_{\theta}
m}\bigr)
S_{\theta} S_{\theta}^T,
\]
where we have used the shorthand
\[
\overline{\partial_{\theta} m} = \int_{\mathcal{X}}\partial
_{\theta} m(x,\theta) h(x) \,\mathrm{d}\mu(x)
\]
along with the observation that
\[
\int\partial_y \phi(e;\theta) \,\mathrm{d}\nu(e) = \int\partial
_{\theta}
\phi(e;\theta) \,\mathrm{d}\nu(e) = \int\partial^2_y \phi(e;
\theta) \,\mathrm{d}\nu(e) = \int\partial^2_{y\theta} \phi
(e;\theta) \,\mathrm{d}
\nu(e) = 0
\]
and some cancelation. We have retained $\phi$ instead of $\phi^*$ in
terms involving $\partial_{\theta}$ for the sake of notational compactness.

We now have that
\begin{eqnarray*}
&&  \int U_J(e,x)^{TT} \phi(e;\theta) h(x) \,\mathrm{d}\nu(e) \,
\mathrm{d}\mu(x)
\\
&&\quad  = (\overline{\partial_{\theta} m} S_{\theta
}
)^{TT} + \int\frac{\partial_{\theta} \phi(e;\theta
)^{TT}}{\phi(e;\theta)} \,\mathrm{d}\nu(e)
\\
&&\qquad{} - \overline{\partial_{\theta} m} S_{\theta} \int
\frac{ \partial_{y} \phi(e;\theta)\, \partial_{\theta} \phi
(e;\theta)^T}{\phi(e;\theta)} \,\mathrm{d}\nu(e) - \int\frac{
\partial
_{\theta} \phi(e;\theta)\, \partial_{y} \phi(e;\theta)^T}{\phi
(e;\theta)} \,\mathrm{d}\nu(e) S_{\theta}^T
\overline{\partial_{\theta} m}^T
\\
&&\qquad{} + \int\!\!\!\int\bigl(\partial_{\theta}m(x;\theta) - \overline{
\partial_{\theta} m}\bigr) e S_{\theta} S_{\theta}^T
S_{\theta} S_{\theta}^T \varepsilon^T \bigl(
\partial_{\theta}m(x;\theta) - \overline{\partial_{\theta}
m}\bigr)
\phi(e;\theta) h(x) \,\mathrm{d}\nu(e) \,\mathrm{d}\mu(x)
\end{eqnarray*}
by making a change of variables $\varepsilon= S_{\theta}^{-1} e$ in the
last line and some cancelation we have that
\begin{eqnarray*}
&& \int U_J(e,x)^{TT} \phi(e;\theta) h(x) \,\mathrm{d}\nu(e) \,
\mathrm{d}\mu(x)
\\
&&\quad  = \int\!\!\!\int\partial_{\theta} m(x;\theta)
\frac
{\partial_y \phi(e;\theta)^{TT}}{\phi(e;\theta)} \partial_{\theta
} m(x;\theta)^T \,\mathrm{d}\mu(x) \,\mathrm{d}
\nu(y) + \int\frac{\partial_{\theta}
\phi(e;\theta)^{TT}}{\phi(e;\theta)} \,\mathrm{d}\nu(e)
\\
&&\qquad{} - \overline{\partial_{\theta} m} \int\frac{
\partial_{y} \phi(e;\theta)\, \partial_{\theta} \phi(e;\theta
)^T}{\phi(e;\theta)} \,\mathrm{d}
\nu(e) - \int\frac{ \partial_{\theta} \phi
(e;\theta)\, \partial_{y} \phi(e;\theta)^T}{\phi(e;\theta)} \,
\mathrm{d}\nu(e) \overline{\partial_{\theta}
m}^T
\end{eqnarray*}
which is readily verified to be the Fisher information for this model.

Where $\theta= (\theta_1,\theta_2)$ can be partitioned into
parameters $\theta_1$ that appear only in $m(x;\theta_1)$ and
parameters $\theta_2$ that appear only in $\phi(e;\theta_2)$ the
terms on the second line above are zero and the resulting information
matrix is diagonal. In the classical case of nonlinear regression with
homoscedastic normal errors, we have
\[
y_i = m(x_i,\theta) + \varepsilon_i,\qquad
\varepsilon_i \sim N\bigl(0,\sigma^2\bigr)
\]
the score covariance for $(\theta,\sigma)$ reduces to
\[
\int U_J(e,x)^{TT} \phi(e;\theta) h(x) \,\mathrm{d}\nu(e) \,\mathrm
{d}\mu(x) =
\lleft[
\matrix{\displaystyle\frac{1}{\sigma^2} \int
\nabla_{\theta} m(x;\theta)^{TT} h(x) \,\mathrm{d}x & 0
\vspace*{3pt}\cr
0 & \displaystyle\frac{1}{2 \sigma^4 }}
\rright].
\]

\textit{Joint centering and conditioning}: $x^s= x$. Here we center and
condition on the entire set of $x$. In this case our results are those
of the uncentered case:
\begin{eqnarray*}
U_J(e,x) & =& J\bigl(e + m(x),x\bigr),
\\
B_n & =& 2\int\!\!\!\int J\bigl(e + m(x),x\bigr) E \hatgnh(x,e,m) \,
\mathrm{d}\nu(e) \,\mathrm{d}
\mu(x).
\end{eqnarray*}
For (\ref{DCLT}):
\begin{eqnarray*}
\tilde{U}_J(e,x) & =& U_J(e,x) + 2 \int J
\bigl(e_*+m(x),x\bigr) f(e_*\mid x) \,\mathrm{d}\nu(e_*),
\\
\tilde{B}_n & =& 2\int\!\!\!\int J\bigl(e + m(x),x\bigr) \frac{E \hatgnh
(x,e)}{E\hathnh(x)}
h(x) \,\mathrm{d}\nu(e) \,\mathrm{d}\mu(x).
\end{eqnarray*}
In this case, $U_J(x,y)$ and $\tilde{U}_J(x,y)$ are exactly the same
as the non-centered case, yielding the information matrix with the
replacement (\ref{eqscorereplace}).

We note that while the non-centered and the jointly centered and
conditioned cases always yield the Fisher information under the
substitution (\ref{eqscorereplace}), the case of centering by some
variables and conditioning on others need not. Even in the
homoscedastic case, efficiency is only gained when the variance of the
model for the residuals is equal to the Fisher information for its
mean. However, under these conditions, we can gain efficiency while
reducing the bias in the central limit theorem above.

Employing these lemmas, we can demonstrate a central limit theorem for
minimum conditional disparity estimates:

%
\begin{teo} \label{generalefficiency}
Let $\{(X_{n1},X_{n2},Y_{n1},Y_{n2}), n \geq1\}$ be given as in
Section~\ref{secassumptions}, under assumptions \textup{(D1)--(D4)},
\textup{(E1)--(E2)}, \textup{(K1)--(K6)},
\textup{(B1)--(B4)}, \textup{(P1)}--\textup{(P4)} and
\textup{(N1)--(N4)} define
\[
\theta_f = \argmin_{\theta\in\Theta} D_\infty(f,\theta)
\]
and
\begin{eqnarray*}
H^D(\theta) & =& \nabla_\theta^2
D_\infty(f,\theta)
\\
& =& \int A_2 \biggl( \frac{f(y\mid x)}{\phi(y\mid x,\theta)} \biggr)
\nabla^2_\theta\phi(y\mid x,\theta) h(x) \,\mathrm{d}\mu(x) \,
\mathrm{d}\nu(y)
\\
&&{} + \int A_3 \biggl( \frac{f(y\mid x)}{\phi(y\mid x,\theta)} \biggr)
\frac{\nabla_\theta
\phi(y\mid x,\theta)^{TT}}{\phi(y\mid x,\theta)} h(x) \,\mathrm
{d}\mu(x) \,\mathrm{d}\nu(y),
\\
I^D(\theta) & =& H^D(\theta) V^D(
\theta)^{-1} H^D(\theta),
\\
\tilde{I}^D(\theta) & =& H^D(\theta)
\tilde{V}^D(\theta)^{-1} H^D(\theta)
\end{eqnarray*}
then
\[
\sqrt{n} \bigl[ T_n(\breve{f}_n) - \theta_f
- B_n \bigr] \rightarrow N \bigl(0,I^D(
\theta_f)^{-1} \bigr)
\]
and
\[
\sqrt{n} \bigl[ \tilde{T}_n(\breve{f}_n) -
\theta_f - \tilde{B}_n \bigr] \rightarrow N \bigl(0,
\tilde{I}^D(\theta_f)^{-1} \bigr)
\]
in distribution where $B_n$, $\tilde{B}_n$, $V^D(\theta)$ and $\tilde
{V}^D(\theta)$
are obtained by substituting
%
\begin{equation}
\label{Jsubs} J(y,x) = A_1 \biggl(\frac{f(y\mid x)}{\phi(y\mid
x,\theta_f)} \biggr)
\frac{\nabla_\theta
\phi(y\mid x,\theta_f)}{\phi(y\mid x,\theta_f)}
\end{equation}
into the expressions for $B_n$, $\tilde{B}_n$, $V_J$ and $\tilde
{V}_J$ in Lemma~\ref{hatfCLT}.
\end{teo}

Here we note that in the case that $f = \phi_{\theta_0}$ for some
$\theta_0$,
that $\theta_f = \theta_0$ and further since $A_1(1)=A_2(1)=A_3(1) =
1$ we have
that $H^D(\theta_f)$ is given by the Fisher information for $\phi
_{\theta_0}$. Since we have demonstrated above that $V^D(\theta_f)$
and $\tilde{V}^D(\theta_f)$ also correspond to the Fisher information
in particular cases above, when this holds $I^D(\theta_f)$ and
$I^{\tilde{D}}(\theta_f)$ also give us the Fisher information and
hence efficiency.

\begin{pf*}{Proof of Theorem \ref{generalefficiency}}
We will define $\bar{T}_n$, and
$f_K(y\mid x)$ to be either the pair ($T_n$, $E [\hatgnh(x,y)]/ \hathnh
(x)$) or ($\tilde{T}_n$, $E
\hatgnh(x,y)/E \hathnh(x)$).
Our arguments now follow those in Tamura and Boos \cite{TamuraBoos86} and Park and Basu
\cite{ParkBasu04}.

Since $\bar{T}_n(f)$ satisfies
\[
\nabla_\theta D_n\bigl(f,\bar{T}_n(f)\bigr) = 0
\]
we can write
\[
\sqrt{n} \bigl( \bar{T}_n(\breve{f}_n) -
\theta_0 \bigr) = - \bigl[ \nabla^2_\theta
D_n\bigl(\breve{f}_n,\theta^+\bigr) \bigr]^{-1}
\sqrt{n} \nabla_\theta D_n(\breve{f}_n,
\theta_0)
\]
for some $\theta^+$ between $\bar{T}_n(\breve{f}_n)$ and $\theta
_f$. It is
therefore sufficient to demonstrate:
\begin{enumerate}[(ii)]
\item[(i)] $\nabla^2_\theta D_n(\breve{f}_n,\theta^+) \rightarrow
H^D(\theta_f)$
in probability.

\item[(ii)]$\sqrt{n} [ \nabla_\theta
D_n(\breve{f}_n,\theta_f) - \bar{B}_n ] \rightarrow
N(0,V^D(\theta_f)^{-1})$ in distribution
\end{enumerate}
with $\bar{B}_n$ given by $B_n$ or $\tilde{B}_n$ as appropriate.

We begin with \textup{(i)} where we observe that by assumption
\textup{(N4)}, $A_2(r)$ and $A_3(r)$ are bounded and the
result follows
from Theorems~\ref{conditionalL1} and~\ref{consistency} and the dominated
convergence theorem. In the case of Hellinger distance
\begin{eqnarray*}
\nabla^2_\theta D(\breve{f}_n,\phi\mid x,
\theta) & =& \int\biggl[ \frac{\nabla^2_\theta
\phi(y,x,\theta)}{\sqrt{\phi(y,x,\theta)}} - \frac{\nabla_\theta
\phi(y,x,\theta)^{TT}}{\phi(y,x,\theta)^{3/2}} \biggr] \sqrt{
\breve{f}_n(y\mid x)}\,\mathrm{d}\nu(y)
\\
& =& \int\sqrt{\phi(y,x,\theta)} \nabla_{\theta} \Psi_\theta(y,x,
\theta) \sqrt{\breve{f}_n(y\mid x)}\,\mathrm{d}\nu(y)
\end{eqnarray*}
so that $| \nabla^2_\theta D_n(\breve{f}_n,\phi\mid x,\theta^+) -
H^D(\theta_f)|$ can be expressed as
\begin{eqnarray*}
&& \int\!\!\!\int\sqrt{\phi(y,x,\theta)} \nabla_\theta\psi(y,x,\theta
) \bigl(
\sqrt{\breve{f}_n(y\mid x)} - \sqrt{f(y\mid x)} \bigr) \,\mathrm
{d}\nu(y)h(x) \,\mathrm{d}
\mu(x)
\\
&&\qquad{} + \int\biggl( \frac{\nabla^2_\theta\phi(y,x,\theta
^+)}{\sqrt{\phi(y,x,\theta^+)}} - \frac{\nabla^2_\theta\phi
(y,x,\theta_f)}{\sqrt{\phi(y,x,\theta_f)}} \biggr)
\sqrt{f(y\mid x)} \,\mathrm{d}\nu(y) h(x) \,\mathrm{d}\mu(x)
\\
&&\quad  \leq\sup_{x \in\mathcal{X}} S \biggl(\int\bigl\vert
\breve{f}_n(y\mid x) - f(y\mid x) \bigr\vert\,\mathrm{d}\nu(y)
\biggr)^{1/2} + \mathrm{o}_p(1)
\\
&&\quad  = \mathrm{o}_p(1).
\end{eqnarray*}
Where the calculations above follow from assumption (N3), bounding (squared) Hellinger distance by $L_1$
distance, the uniform $L_1$ convergence of $\breve{f}_n$
(Theorem~\ref{jointsupL1}) and the consistency of $\theta$ (Theorem
\ref{consistency}).

Turning to \textup{(ii)} where we observe that by the boundedness
of $C$
and the dominated convergence theorem, we can write $\nabla_\theta
D_n(\breve{f}_n,\phi\mid x,\theta) - \bar{B}_{n}$ as
\begin{eqnarray*}
&& \int A_2 \biggl( \frac{\breve{f}_n(y\mid x)}{\phi(y\mid x,\theta
)} \biggr) \nabla_\theta
\phi(y\mid x,\theta) \,\mathrm{d}\nu(y) - \bar{B}_{n}
\\
&&\quad  = \int A_1 \biggl( \frac{f_K(y\mid x)}{\phi(y\mid x,\theta)}
\biggr)
\frac{\nabla_\theta\phi(y\mid x,\theta)}{\phi(y\mid x,\theta)}
\bigl[ \breve{f}_n(y\mid x) - f_K(y
\mid x) \bigr]\,\mathrm{d}\nu(y)
\\
&&\qquad{}  + \int\biggl[ A_2 \biggl( \frac{\breve{f}_n(y\mid x)}{\phi
(y\mid x,\theta)} \biggr)
- A_2 \biggl( \frac{f_K(y\mid x)}{\phi(y\mid x,\theta)} \biggr)
\biggr]
\nabla_\theta\phi(y\mid x,\theta) \,\mathrm{d}\nu(y)
\\
&&\qquad{} - \int A_1 \biggl( \frac{f_K(y\mid x)}{\phi(y\mid x,\theta)}
\biggr)
\biggl( \frac{
\breve{f}_n(y\mid x)}{\phi(y\mid x,\theta)} - \frac{f_K(y\mid
x)}{\phi(y\mid x,\theta)} \biggr) \nabla_\theta
\phi(y\mid x,\theta) \,\mathrm{d}\nu(y)
\end{eqnarray*}
from a minor modification Lemma~25 of Lindsay \cite{Lindsay94} we have
that by the
boundedness of $A_1$ and $A_2$ there is a constant $B$ such that
\[
\bigl\vert A_2\bigl(r^2\bigr) - A_2
\bigl(s^2\bigr) - \bigl(r^2-s^2
\bigr)A_1\bigl(s^2\bigr)\bigr\vert\leq
\bigl(r^2 - s^2\bigr)B
\]
substituting
\[
r = \sqrt{\frac{\breve{f}_n(y\mid x)}{\phi(y\mid x,\theta)}},\qquad s=
\sqrt{\frac{f_K(y\mid x)}{\phi(y\mid x,\theta)}}
\]
we obtain
\begin{eqnarray*}
&&  \int A_2 \biggl( \frac{\breve{f}_n(y\mid x)}{\phi(y\mid x,\theta
)} \biggr) \nabla_\theta
\phi(y\mid x,\theta) \,\mathrm{d}\nu(y) - \bar{B}_{n}
\\
&&\quad  = \int A_1 \biggl( \frac{f_K(y\mid x)}{\phi(y\mid x,\theta)}
\biggr)
\frac{\nabla_\theta\phi(y\mid x,\theta)}{\phi(y\mid x,\theta)}
\bigl[ \breve{f}_n(y\mid x) - f_K(y
\mid x) \bigr]\,\mathrm{d}\nu(y)
\\
&&\qquad{}  + B\int\frac{\nabla_\theta
\phi(y\mid x,\theta)}{\phi(y\mid x,\theta)} \bigl( \sqrt{\breve{f}_n(y
\mid x)} - \sqrt{f_K(y\mid x)} \bigr)^2 \,\mathrm{d}\nu(y).
\end{eqnarray*}
The result now follows from Lemmas~\ref{hatfHD} and~\ref{hatfCLT}.

For the special case of Hellinger distance, we observe that
\[
\nabla_\theta D_n(\breve{f}_n,\phi\mid x,
\theta) = \int\frac{\nabla
_\theta
\phi(y,x,\theta)}{\sqrt{\phi(y,x,\theta)}}\sqrt{\breve
{f}_n(y\mid x)}\,\mathrm{d}\nu(y)
\]
and\vspace*{2pt} applying the identity $\sqrt{a}-\sqrt{b} =
(a-b)/2\sqrt{a} + (\sqrt{b}-\sqrt{a})^2/2\sqrt{a}$ with $a =
f_K(y\mid x)$ and $b = \breve{f}_n(y\mid x)$, we obtain
\begin{eqnarray*}
&& \sqrt{n} \int\frac{\nabla_\theta
\phi(y,x,\theta)}{\sqrt{\phi(y,x,\theta)}} \bigl( \sqrt{\breve
{f}_n(y\mid
x)} - \sqrt{f_K(y\mid x)} \bigr) \,\mathrm{d}\nu(y)
\\
&&\quad  = \sqrt{n} \int\frac{\nabla_\theta
\phi(y,x,\theta)}{2 \sqrt{\phi(y,x,\theta) f_K(y\mid x)}} \bigl
(\breve
{f}_n(y\mid x) - f_K(y\mid x) \bigr) \,\mathrm{d}\nu(y)
\\
&&\qquad{} - \sqrt{n} \int\frac{\nabla_\theta
\phi(y,x,\theta)}{2 \sqrt{\phi(y,x,\theta) f_K(y\mid x)}} \bigl(
\sqrt{
\breve{f}_n(y\mid x)} - \sqrt{f_K(y\mid x)}
\bigr)^2 \,\mathrm{d}\nu(y)
\\
&&\quad  = \sqrt{n} \biggl( \int\frac{\nabla_\theta
\phi(y,x,\theta)}{2 \sqrt{\phi(y,x,\theta) f_K(y\mid x)}} \breve
{f}_n(y\mid x) - B_n \biggr) + \mathrm{o}_p(1),
\end{eqnarray*}
where we have applied Lemma~\ref{hatfHD} to the second term in the
expression above, and can now obtain the result from Lemma~\ref
{hatfCLT} and the convergence of $f_K(y\mid x)$ to $f(y\mid x)$.
\end{pf*}

We note here that Theorem~\ref{generalefficiency} relies on
assumption \textup{(D4)} only through the consistency of $\bar
{T}_n(\breve{f}_n)$ and Lemmas~\ref{hatfHD} and~\ref{hatfCLT}. In
the case of $\tilde{T}_n(\hat{f}_n^{*})$ (uncentered densities with
the integral form of the disparity), we can remove this condition by
employing Theorem~E.1, and Lemmas~E.1
and~E.2 from Supplemental Appendix~E
(Hooker \cite{Hooker14}).

\section{Robustness properties} \label{secrobustness}

An important motivator for the study of disparity methods is that in addition
to providing statistical efficiency as demonstrated above, they are
also robust
to contamination from outlying observations. Here we investigate the robustness
of our estimates through their breakdown points.
These have been studied for i.i.d. data in Beran \cite{Beran77};
Park and Basu \cite{ParkBasu04};
Lindsay \cite{Lindsay94}
and the extension to conditional models follows similar lines.

In particular, we examine two models for contamination:
\begin{enumerate}[3.]
\item To mimic the ``homoscedastic'' case, we contaminate
$g(x_1,x_2,y_1,y_2)$ with
outliers independent of $(x_1,x_2)$. That is, we define the contaminating
density
%
\begin{equation}
\label{contamination1} g_{\varepsilon,z}(x_1,x_2,y_1,y_2)
= (1-\varepsilon) g(x_1,x_2,y_1,y_2)
+ \varepsilon\delta_z(y_1,y_2)
h(x_1,x_2),
\end{equation}
where $\delta_z$ is a contamination
density parameterized by $z$ such that $\delta_z$ becomes ``outlying''
as $z \rightarrow\infty$. Typically, we think of $\delta_z$ as
having small
support centered around $z$. This results in the conditional density
\[
f_{\varepsilon,z}(y_1,y_2\mid x_1,x_2)
= (1-\varepsilon)f(y_1,y_2\mid x_1,x_2)
+ \varepsilon\delta_z(y_1,y_2)
\]
which we think of as the result of smoothing a contaminated residual density.
We note that we have not changed the marginal distribution of
$(x_1,x_2)$ via
this contamination. This particularly applies to the case where only
$y_1$ is present and
the estimate (\ref{hatm})--(\ref{hatft}) is employed.

\item In the more general setting, we set
%
\begin{equation}
\label{contamination2} g_{\varepsilon,z}(x_1,x_2,y_1,y_2)
= (1-\varepsilon) g(x_1,x_2,y_1,y_2)
+ \varepsilon\delta_z(y_1,y_2)
J_U(x_1,x_2) h(x_1,x_2),
\end{equation}
where $J_U(x_1,x_2)$ is the indicator of $(x_1,x_2) \in U$ scaled so that
$h(x_1,x_2)J_U(x_1,x_2)$ is a distribution. This translates to the conditional
density
\[
f_{\varepsilon,z}(y_1,y_2\mid x_1,x_2)
= \cases{ f(y_1,y_2\mid x_1,x_2),
&\quad$(x_1,x_2) \notin U$,
\cr
(1-
\varepsilon)f(y_1,y_2\mid x_1,x_2)
+ \varepsilon\delta_z(y_1,y_2), &\quad
$(x_1,x_2) \in U$}
\]
which localizes contamination in covariate space. Note that the marginal
distribution is now scaled differently in $U$.
\end{enumerate}
Naturally, this characterization (\ref{contamination1}) does not
account for the effect of outliers on
the Nadaraya--Watson estimator (\ref{hatm}). If these are localized in
covariate
space, however, we can think of (\ref{hatft}) as being approximately a mixture
of the two cases above. As we will see the distinction between these
two will
not affect the basic properties below. Throughout we will write
$\delta_z(y_1,y_2\mid x_1,x_2)$ in place of $\delta_z(y_1,y_2)$ or
$\delta_z(y_1,y_2)J_U(x_1,x_2)$ as appropriate. $h(x_1,x_2)$ will be
taken to be
modified according to (\ref{contamination2}) if appropriate.

We must first place some conditions on $\delta_z$:
\begin{enumerate}[C1.]
\item[C1.]$\delta_z$ is orthogonal in the limit to $f$. That is
\[
\lim_{z \rightarrow\infty} \sum_{y_2 \in S_y} \int\delta
_z(y_1,y_2\mid x_1,x_2)
f(y_1,y_2\mid x_1,x_2)
\,\mathrm{d}y_1 = 0\qquad \forall(x_1,x_2).
\]

\item[C2.]$\delta_z$ is orthogonal in the limit to $\phi$:
\[
\lim_{z \rightarrow\infty} \sum_{y_2 \in S_y} \int\delta
_z(y_1,y_2\mid x_1,x_2)
\phi(y_1,y_2\mid x_1,x_2,
\theta) \,\mathrm{d}y_1 = 0\qquad \forall(x_1,x_2).
\]

\item[C3.]$\phi$ becomes orthogonal to $f$ for large $\theta$:
\[
\lim_{\llVert \theta\rrVert \rightarrow\infty} \sum_{y_2 \in
S_y} \int
f(y_1,y_2\mid x_1,x_2)
\phi(y_1,y_2\mid x_1,x_2,\theta)
= 0\qquad \forall(x_1,x_2).
\]

\item[C4.]$C(-1)$ and $C'(\infty)$ are both finite or the disparity is Hellinger
distance.
\end{enumerate}

In the following result with use $T[f] = \argmin D_\infty(f,\theta)$
for any
$f$ in place of our estimate~$\hat{\theta}$.

%
\begin{teo} \label{RobustnessTheorem1}
Under assumptions \textup{C1--C4}
under both
contamination models (\ref{contamination1}) and (\ref
{contamination2}) define
$\varepsilon^*$ to satisfy
%
\begin{equation}
\label{epscondition} \bigl(1-2\varepsilon^*\bigr) C'(\infty) =
\inf
_{\theta\in\Theta} D \bigl( \bigl(1-\varepsilon^*\bigr) f,\theta
\bigr) - \lim
_{z \rightarrow\infty} \inf_{\theta\in
\Theta} D \bigl( \varepsilon^*
\delta_z, \theta\bigr)
\end{equation}
with $C'(\infty)$ replaced by 1 in the case of Hellinger distance then for
$\varepsilon< \varepsilon^*$
\[
\lim_{z \rightarrow\infty} T[f_{\varepsilon,z}] = T\bigl
[(1-\varepsilon)f\bigr]
\]
and in particular the breakdown point is at least $\varepsilon^*$: for
$\varepsilon< \varepsilon^*$,
\[
\sup_z \bigl\llVert T[f_{\varepsilon,z}] - T\bigl[(1-
\varepsilon) f\bigr] \bigr\rrVert< \infty.
\]
\end{teo}


\begin{pf}
We begin by observing that by assumption C1, for any
fixed $\theta$,
\begin{eqnarray*}
D( f_{\varepsilon,z},\theta) & =& \int\!\!\!\int_{A_z(x)} C \biggl(
\frac{f_{\varepsilon,z}(y\mid x)}{\phi(y\mid x,\theta)} - 1 \biggr
) \phi(y\mid x,\theta) h(x) \,\mathrm{d}\nu(y) \,\mathrm{d}\mu(x)
\\
&&{} + \int\!\!\!\int_{A_z^c(x)} C \biggl( \frac{f_{\varepsilon,z}(y\mid
x)}{\phi(y\mid x,\theta)}
- 1 \biggr) \phi(y\mid x,\theta) h(x) \,\mathrm{d}\nu(y) \,\mathrm
{d}\mu(x)
\\
& =& D_{A_z}(f_{\varepsilon,z},\theta) + D_{A_z^c}(f_{\varepsilon,z},
\theta),
\end{eqnarray*}
where $A_z(x) = \{y\dvt \max(f(y\mid x),\phi(y\mid x,\theta)) >
\delta_z(y\mid x) \}$. We note that for any $\eta$ with $z$ sufficiently
large that
\[
\sup_{x \in\mathcal{X}} \sup_{y \in A_z(x)}
\delta_z(y\mid x) < \eta\quad\mbox{and}\quad \sup_{(x) \in\mathcal{X}}
\sup_{y \in A_z^c(x)} f(y\mid x) < \eta
\]
and thus for sufficiently large $z$,
\begin{eqnarray*}
&& \bigl\vert D( f_{\varepsilon,z},\theta) - \bigl(D_{A_z}\bigl((1-
\varepsilon)f,\theta\bigr) + D_{A_z^c}(\varepsilon\delta_z,\theta)
\bigr) \bigr\vert
\\
&&\quad  \leq \int\!\!\!\int C \biggl(\frac{\eta}{\phi
(y\mid x,\theta)}-1 \biggr) \phi(y
\mid x,\theta)h(x) \,\mathrm{d}\nu(y) \,\mathrm{d}\mu(x)
\\
&&\quad  \leq\eta\sup_t \bigl\vert C'(t)\bigr
\vert
\end{eqnarray*}
hence
%
\begin{equation}
\label{Dbreakup} \sup_{\theta} \bigl\vert D( f_{\varepsilon,z},
\theta) - \bigl(D_{A_z}\bigl((1-\varepsilon)f,\theta\bigr) +
D_{A_z^c}(\varepsilon\delta_z,\theta) \bigr) \bigr\vert
\rightarrow0.
\end{equation}
We also observe that for any fixed $\theta$,
\begin{eqnarray*}
D_{A_z^c}(\varepsilon\delta_z,\theta) & =& \int\!\!\!\int
_{A_z^c(x)} C \biggl(\frac{2\varepsilon\delta_z(y)}{\phi
(y\mid x,\theta)}-1 \biggr) \phi(y\mid x,
\theta) h(x) \,\mathrm{d}\nu(y) \,\mathrm{d}\mu(x)
\\
&&{} + \int\!\!\!\int_{A_z^c(x)} \varepsilon
\delta_z(y\mid x) C'\bigl( t(y,x) \bigr) \,\mathrm{d}\nu(y) \,
\mathrm{d}
\mu(x)
\\
& \rightarrow&\varepsilon C'(\infty)
\end{eqnarray*}
for $t(y,x)$ between $\varepsilon
\delta_z(y\mid x)/\phi(y\mid x,\theta)$ and $2\varepsilon
\delta_z(y\mid x)/\phi(y\mid x,\theta)$ since $t(y,x) \rightarrow
\infty$,
$C(\cdot)$ and $C'(\cdot)$ are bounded
and $\int_{A_z^c(x)} \phi(y\mid x,\theta) \,\mathrm{d}\nu(y)
\rightarrow0$.

Similarly,
\[
D_{A_z^c}\bigl((1-\varepsilon)f,\theta\bigr) \rightarrow D\bigl((1-
\varepsilon)f,\theta\bigr)
\]
and thus
\[
D(f_{\varepsilon,z},\theta) \rightarrow D\bigl((1-\varepsilon
)f,\theta\bigr) +
\varepsilon C'(\infty)
\]
which is minimized at $\theta= T[f_{\varepsilon,z}]$.

It remains to rule out divergent sequences $\llVert \theta_z\rrVert
\rightarrow
\infty$.
In this case, we define $B_z(x) = \{y\dvt f(y\mid x) > \max
(\varepsilon
\delta_z(y\mid x), \phi(y\mid x,\theta_z)) \}$ and note that from the
arguments above
\[
D_{B_z}\bigl((1-\varepsilon)f,\theta_z\bigr) \rightarrow(1-
\varepsilon) C'(\infty)
\]
and
\[
D_{B_z^c}(\varepsilon\delta,\theta_z) \rightarrow D(\varepsilon
\delta,\theta_z)
\]
and hence
\[
\lim_{z \rightarrow\infty} D(f_{\varepsilon,z},\theta_z) > \lim
_{z
\rightarrow
\infty} \inf_{\theta\in\Theta} D(\varepsilon
\delta_z,\theta) + (1-\varepsilon) C'(\infty) > D
\bigl(f_{\varepsilon,z},T\bigl[(1-\varepsilon) f\bigr] \bigr)
\]
from (\ref{epscondition}), yielding a contradiction.

In the case of Hellinger distance, we observe
\begin{eqnarray*}
&& \bigl\vert D( f_{\varepsilon,z},\theta) - \bigl(D\bigl
((1-\varepsilon)f,\theta
\bigr) + D(\varepsilon\delta_z,\theta) \bigr) \bigr\vert
\\
&&\quad  = \int\!\!\!\int\sqrt{\phi(y\mid x,\theta)} \bigl(
\sqrt{f_{\varepsilon,z}(y\mid x)} - \sqrt{(1-\varepsilon)f(y\mid
x)} - \sqrt{\varepsilon
\delta_z(y\mid x)} \bigr) h(x) \,\mathrm{d}\nu(y) \,\mathrm{d}\mu(x)
\\
&&\quad  \leq\int\!\!\!\int\bigl( \sqrt{f_{\varepsilon,z}(y\mid x)} - \sqrt{(1-
\varepsilon)f(y\mid x)} - \sqrt{\varepsilon\delta_z (y\mid x)}
\bigr)^2 h(x) \,\mathrm{d}\nu(y) \,\mathrm{d}\mu(x)
\\
&&\quad  = \int\!\!\!\int\bigl[ 2(1-\varepsilon)f(y\mid x) + 2 \varepsilon
\delta(y
\mid x) \bigr] h(x) \,\mathrm{d}\nu(y) \,\mathrm{d}\mu(x)
\\
&&\qquad{} - 2\int\!\!\!\int\bigl( \sqrt{f_{\varepsilon,z}(y\mid x)} \bigl(
\sqrt{(1-\varepsilon) f(y\mid x)} + \sqrt{\varepsilon\delta
_z(y\mid x)}
\bigr) \bigr) h(x) \,\mathrm{d}\nu(y) \,\mathrm{d}\mu(x),
\end{eqnarray*}
where, by dividing the range of $y$ into $A_z(x)$ and $A_z^c(x)$
as above, we find that on $A_z(x)$, for any $\eta> 0$ and $z$
sufficiently large,
\begin{eqnarray*}
&& \bigl\vert(1-\varepsilon)f(y\mid x) - \sqrt{f_{\varepsilon
,z}(y\mid x)} \bigl(
\sqrt{(1-\varepsilon) f(y\mid x)} + \sqrt{\varepsilon\delta
_z(y\mid x) }
\bigr)\bigr\vert\leq\sqrt{\varepsilon\eta f_{\varepsilon
,z}(y\mid x)} + \varepsilon
\eta
\end{eqnarray*}
which with the corresponding arguments on $A_z^c(x)$ yields
(\ref{Dbreakup}). We further observe that for fixed~$\theta$
\[
D(\varepsilon\delta_z,\theta) = 1 + \varepsilon- \sqrt{\varepsilon
} \int
\sqrt{\delta_z(y\mid x) \phi(y\mid x,\theta)} h(x) \,\mathrm{d}\nu
(y) \,\mathrm{d}\mu(x)
\rightarrow1 + \varepsilon
\]
and for $\llVert \theta_z\rrVert \rightarrow\infty$,
\[
D\bigl((1-\varepsilon)f,\theta_z\bigr) = 2- \varepsilon- \sqrt{1-
\varepsilon} \int\sqrt{f(y\mid x) \phi(y\mid x,\theta_z)} h(x) \,
\mathrm{d}\nu(y) \,\mathrm{d}
\mu(x) \rightarrow2- \varepsilon
\]
from which the result follows from the same arguments as above.
\end{pf}

These results extend on Park and Basu  \cite{ParkBasu04} and Beran \cite{Beran77} and a number of
ways and a few remarks are warranted:
\begin{enumerate}[2.]
\item In Beran \cite{Beran77}, $\Theta$ is assumed to be compact, allowing
$\theta_z$ to converge at least on a subsequence. This removes the
$\llVert \theta_z\rrVert \rightarrow\infty$ case and the result
can be shown for
$\varepsilon\in[0, 1)$.

\item We have not assumed that the uncontaminated density $f$ is a
member of the parametric class
$\phi_\theta$. If $f = \phi_{\theta_0}$ for some $\theta_0$, then
we observe
that by Jensen's inequality
\[
D\bigl((1-\varepsilon)\phi_{\theta_0},\theta\bigr) >
C(-\varepsilon) = D\bigl((1-
\varepsilon)\phi_{\theta_0},\theta_0\bigr)
\]
hence $T[(1-\varepsilon) \phi_{\theta_0}] = \theta_0$. We can
further bound
$D(\varepsilon\delta_z, \theta) > C(\varepsilon-1)$ in which case
(\ref{epscondition}) can be bounded by
\[
(1-2\varepsilon) C'(\infty) \geq C(-\varepsilon) - C(\varepsilon-1)
\]
which is satisfied for $\varepsilon= 1/2$. We note that in the more general
condition, if $(1-\varepsilon)f$ is closer to the family $\phi_\theta
$ than
$\varepsilon\delta_z$ at $\varepsilon= 1/2$, the breakdown point
will be greater
than $1/2$; in the reverse situation it will be smaller.
\end{enumerate}

We emphasize here that we consider robustness here in the sense of
having outliers in the response variables $Y_i$. Outliers in the $X_i$
result in points of high leverage, to which our methods are not robust.
Robustness in this sense would require a weighted combination of the
$D_n(f,\phi\mid x,\theta)$ as an objective and the resulting efficiency
properties of the model are not clear.


\section{Bandwidth selection, bootstrapping, bias correction and inference} \label{secbootstrap}

The results in the previous sections indicate that minimum disparity
estimates based on non-parametric conditional density estimates are
efficient in the sense that their asymptotic variance is identical to
the Fisher information when the model is correct. They are also robust
to outliers. This comes at a price, however, of a bias that is
asymptotically non-negligible. Here, we propose to correct this bias
with a bootstrap based on the estimated conditional densities. This
will also provide a means of inference that does not assume the
parametric model. We also provide details of the bandwidth selection
methods used in our empirical studies. The details in this section are
heuristic choices applied to the simulation studies in Section~\ref
{secsimulations} and real data analysis in Section~\ref{secdata}.

\subsection{Bandwidth selection}

Bandwidth selection is not particularly well studied for multivariate
or conditional density estimates and software implementing existing
methods is not readily available. Here, we employed a na\"{i}ve
cross-validation approach designed to be methodologically
straightforward. In particular:
\begin{enumerate}[3.]
\item We chose bandwidths $\cnxmb$ for $\hatmnh$ by cross-validating
squared error.

\item We chose bandwidths $\cnxgb$ associated with $x^{\bar{g}}$ in
$\hathnh$ by cross-validating the non-parametric log likelihood:
\[
\cnxgb= \argmax\sum_{i=1}^n \log
\hathnh^{-i}\bigl(X^{\bar{g}}_i\bigr),
\]
where $\hathnh^{-i}$ is the estimate $\hathnh$ based on the data set
with $X^{\bar{g}}_i$ removed.

\item We fixed $\hatmnh$ and $\hathnh$ and their bandwidths and chose
$c_{ny}$ based on cross-validating the non-parametric conditional log
likelihood:
\[
c_{ny} = \argmax\sum_{i=1}^n
\log\hatgnh^{-i}\bigl( Y_i - \hatmnh\bigl(X^{\bar{m}}_i
\bigr), X^{\bar{g}}_i\bigr).
\]
Noting that the denominator in the conditional density becomes an
additive term after taking logs and does not change with $c_{ny}$.
\end{enumerate}
Where we also used discrete values $X_2$, these bandwidths were
estimated for each value of $X_2$ separately at each step. The
resulting bandwidths were then averaged in order to improve the
stability of bandwidth selection.

\subsection{Bootstrapping}

We have two aims in bootstrapping: bias correction and inference.
Nominally, we can base inference on the asymptotic normality results
established in Theorem~\ref{generalefficiency} using the inverse of
the Fisher information as the variance for the estimated parameters.
However, the coverage probabilities of confidence intervals based on
these results will be poor due to the non-negligible bias in the
theorem; it will also not provide correct coverage when the assumed
parametric model is incorrect.

As an alternative, we propose a bootstrap based on the estimated
non-parametric conditional densities. That is, to create each bootstrap
sample, we simulate a new response $Y^*_i$ from $\breve{f}_n(\cdot
|X_{i1},X_{i2})$ for $i = 1,\ldots,n$ and use these to re-estimate
parameters $\hat{\theta}$. For continuous $Y_{i1}$, simulating from
this density can be achieved by choosing $Y_{j1}$ with weights $K(
[X_i- X_j]/\cnxgb)$ and then simulating from the density
$c_{ny}^{-\mathrm{d}y}
K( (y-Y_i)/c_{ny})$. For discrete $Y_{i2}$, simulating from the
non-parametric multinomial model is straightforward.

In the simulation experiments below, we examine a number of different
choices of $x^{\bar{m}}$ and $x^{\bar{g}}$ and each is bootstrapped
separately. For maximum likelihood and other robust estimators, we
employ a residual bootstrap for continuous responses and a parametric
bootstrap for discrete responses.

We also examine a hybrid method proposed in Hooker and Vidyashankar  \cite{HookerVidyashankar13} in which we replace $\hatmnh$ with a parametric
regression model $m(x,\theta)$. We then minimize the disparity between
the estimated density of residuals (which varies with parameters) and a
parametric residual density. Specifically, we set
\begin{eqnarray*}
E_i(\theta) & =& Y_i - m(X_i,\theta),
\\
\tilde{f}_n(e,\theta) & =& \frac{1}{nc_n} \sum
_{i=1}^n K \biggl( \frac{e - E_i(\theta)}{c_n} \biggr),
\\
\tilde{\theta}_n & =& \argmin_{\theta\in\Theta} \int C \biggl(
\frac{\tilde{f}_n(e,\theta)}{\phi(e)} -1 \biggr) \phi(e) \,
\mathrm{d}e.
\end{eqnarray*}
This formulation avoids conditional density estimation (and hence
asymptotic bias) at the expense of a parameter-dependent kernel density
estimate for the residuals. In this formulation $\phi(e)$ is a
reference residual density in which a scale parameter has been robustly
estimated. In the simulations below, the scale parameter is
re-estimated via a disparity method with the remaining $\theta$ held
fixed. For this case, we employ a parametric bootstrap at the estimated
parameters, but sample from the estimated non-parametric residual
density. Throughout, we keep the estimated bandwidths fixed.

\subsection{Inference}

Given a bootstrap sample $\theta_b^*$, $b = 1,\ldots,B$ along with
our original estimate $\hat{\theta}$, we conduct inference along well
established lines:
\begin{itemize}
\item Obtain a bias corrected estimate
\[
\hat{\theta}^c = 2 \hat{\theta} - \frac{1}{B} \sum
_{b=1}^B \theta_b^*.
\]

\item Estimate a bootstrap standard error, $\widehat{\mbox{se}}(\theta)$, from the sample standard deviation of $\theta_b$.

\item Construct confidence intervals $[\hat{\theta}^c - 1.96 \widehat{\mbox{se}}(\theta), \hat{\theta}^c + 1.96 \widehat{\mbox{se}}(\theta)]$.
\end{itemize}
The performance of these confidence intervals will be examined in the
simulation studies below, but we make a couple of remarks on this:
\begin{enumerate}[3.]
\item Our bootstrap scheme amounts to simulation under the model
$\breve{f}_n$. Given the convergence of $\breve{f}_n$ to $f$ in
Theorem~\ref{jointsupL1} and the continuity of $I^D(\theta)$ and
$\tilde{I}^D(\theta)$ in $f$, the bootstrap standard error can be
readily shown to be consistent for the sampling standard error of $\hat
{\theta}$. Similarly, since density estimates with bandwidths $c_{ny}$
and $2c_{ny}$ converge, the bias correction incurs no additional variance.

\item The bias correction for the proposed bootstrap approximates
considering the difference between estimating $\breve{f}_n$ with
bandwidths $c_{ny}$ and $2 c_{ny}$; this is exactly true when employing
a Gaussian kernel. The bias terms in Lemma~\ref{hatfCLT} are readily
shown to be $\mathrm{O}(c_{ny}^2)$ which would suggest a corrected estimate of
the form $(4\hat{\theta} - 1/B \sum\theta_b^*)/3$ instead of the
linear correction proposed above. However\vspace*{1pt} the estimate is also biassed
due to the nonlinear dependence of $\hat{\theta}$ on $\breve{f}_n$
regardless of the value of $c_{ny}$. This bias is asymptotically
negligible, but we have found the proposed correction to provide better
performance at realistic sample sizes. A~combined bias correction
associated with explicitly obtaining an estimate at $2c_{ny}$ to
correct for smoothing bias with a bootstrap estimate to correct for
intrinsic bias may improve performance further, but this is beyond the
scope of this paper.
\end{enumerate}

\section{Simulation studies} \label{secsimulations}

Here we report simulation experiments designed to evaluate the methods
analyzed above. Our examples are all based on conditionally-specified
regression models. In all of these, we generate a three-dimensional set
of covariates in the following manner:
\begin{enumerate}
\item Generate $n\times3$ matrix $X$ from a Uniform random variable on $[-1,1]$.

\item Post-multiply this matrix by a $\sqrt{8}/3$ times a matrix with
unit diagonal and 0.25 in all off-diagonal entries to create correlation.

\item Replace the third column of $X$ with the indicator of the
corresponding entry being greater than zero.
\end{enumerate}
This gives us two continuous valued covariates and a categorical
covariate all of which are correlated. The values of these covariates
were regenerated in each simulation.

Using these covariates, we simulated data from two models:
\begin{itemize}
\item A linear regression with Gaussian errors and all coefficients
equal to 1:
%
\begin{equation}
\label{linregmod} Y_i = 1 + \sum_{j=1}^3
X_{ij} + \varepsilon_i
\end{equation}
with $\varepsilon_i \sim N(0,1)$, This yields a signal to noise ratio of
1.62. In this model, we estimate the intercept and all regression
parameters as well as the noise variance, yielding true values of
$(\beta_0,\beta_1,\beta_2,\beta_3,\sigma) = (1,1,1,1,1)$. We
optimize over $\log\sigma$ to avoid boundary problems and have
reported estimate and standard errors for $\log\sigma$ below.

\item A logistic regression with zero intercept and all other
coefficients 0.5:
%
\begin{equation}
\label{logitmod} P(X_i = 1\mid X_i) =
\frac{\mathrm{e}^{\sum_{j=1}^3 0.5 X_{ij}}}{ 1 + \mathrm{e}^{\sum
_{j=1}^3 0.5 X_{ij}}}
\end{equation}
in order to evaluate a categorical response model. Here only the four
regression parameters were estimated.
\end{itemize}
In each model we also examined the addition of outliers. In (\ref
{linregmod}), we changed either 1, 3, 5 or 10 of the $\varepsilon_i$ to
take values 3, 5, 10 and 15. These covariate values $X_i$ corresponding
the modified $\varepsilon_i$ where held constant within each simulation
study, but were selected in two different ways:
\begin{enumerate}
\item At random from among all the data.

\item Based on the points with $X_{i1}$ closest to $-$0.5.
\end{enumerate}
These mimic the contamination scenarios above.


In binary response data in (\ref{logitmod}), we require a model in
which an ``outlier'' distribution can become orthogonal to the model
distribution. For binary data this can occur only if the parametric
model has $P(Y=1\mid X) \approx0$ or $P(Y=1\mid X) \approx1$ which for
logistic regression can occur only at values of $X$ that have high
leverage; a robustness problem not considered in this paper. Instead,
we examine a logistic binomial model based on successes out of 8
trials. For this, we have employed an exact distribution which is
contaminated with $\alpha\%$ of a distribution in which points take
the value 8, either uniformly as in scenario (\ref{contamination1}) or
at the single $X_i$ with $X_{i1}$ closest to $-$0.5 as in scenario~(\ref{contamination2}). In this case, reasonable estimates of conditional
distributions would require very large sample sizes and we have based
all our estimates on exact distributions.

\subsection{Linear regression}

For the linear regression simulations, we employed 31 points generated
as above. We considered three types of density estimates corresponding
to no centering (labeled HD and NED for Hellinger distance and negative
exponential disparity), jointly centering and conditioning on all
variables (HD.c and NED.c) and the homoscedastic model: centering by
all variables but assuming a constant residual density (HD.h and
NED.h). We also included the marginal method of Hooker and Vidyashankar \cite
{HookerVidyashankar13} which involves only fitting a kernel density
estimate to the residuals of a linear regression. Bandwidths where
chosen by cross-validated log likelihood for uncontaminated data. We
conducted all estimates by minimizing $D_n(\breve{f},\theta)$ with
$D(\breve{f},\phi\mid X_i,\theta)$ approximated a Monte Carlo integral
based on 101 points drawn from $\breve{f}(\cdots\mid X_i)$.

%
\begin{table}[b]
\tabcolsep=0pt
\caption{Simulation results for a linear regression simulation.
\texttt{Lik} are the maximum likelihood estimates, \texttt{G--Y}
correspond to Gervini and Yohai's adaptive truncation estimator,
\texttt{HD} is minimum Hellinger distance, \texttt{NED} is minimum
negative exponential disparity based on uncentered kernel density
estimates, \texttt{HD.c} and \texttt{NED.c} are centered by a
Nadaraya--Watson estimator, \texttt{HD.h} and \texttt{NED.h} are based
on homoscedastic conditional density estimates and \texttt{HD.m} and
\texttt{NED.m} are the marginal estimators in Hooker and Vidyashankar~\cite
{HookerVidyashankar13}. We report the mean value over 5000 simulations
as well as the standard deviation (sd) between simulations} \label{linregtable}
\begin{tabular*}{\tablewidth}{@{\extracolsep{\fill}}@{}llllllllllll@{}} \hline
& $\log\sigma$ & sd & $\beta_0$ & sd & $\beta_1$ & sd & $\beta_2$ & sd & $\beta_3$ & sd & Time \\
\hline
Lik & $-$0.10 & 0.14 & 1.00 & 0.28 & 1.00 & 0.40 & 0.99 & 0.40 & 0.99 & 0.43 & 0.0049\\
G--Y & $-$0.10 & 0.19 & 1.00 & 0.30 & 1.00 & 0.43 & 0.99 & 0.42 & 0.99 & 0.46 & 0.0144\\
HD.c & \phantom{$-$}0.13 & 0.44 & 0.97 & 0.34 & 0.94 & 0.60 & 0.94 & 0.50 & 1.05 & 0.60 & 0.0588\\
NED.c & \phantom{$-$}0.12 & 0.23 & 0.98 & 0.30 & 0.95 & 0.40 & 0.94 & 0.40 & 1.04 & 0.45 & 0.0751\\
HD & \phantom{$-$}0.26 & 0.37 & 0.94 & 0.40 & 0.87 & 0.39 & 0.86 & 0.51 & 1.11 & 0.58 & 0.0604\\
NED & \phantom{$-$}0.25 & 0.16 & 0.94 & 0.30 & 0.87 & 0.35 & 0.87 & 0.36 & 1.11 & 0.45 & 0.0776\\
HD.h & $-$0.18 & 0.32 & 0.95 & 0.50 & 0.88 & 0.34 & 0.88 & 0.34 & 1.10 & 0.75 & 0.0616\\
NED.h & $-$0.17 & 0.21 & 0.95 & 0.34 & 0.88 & 0.34 & 0.88 & 0.34 & 1.09 & 0.50 & 0.0764\\
HD.m & \phantom{$-$}0.05 & 0.17 & 1.00 & 0.29 & 1.00 & 0.43 & 1.00 & 0.42 & 1.00 & 0.45 & 0.0328\\
NED.m & \phantom{$-$}0.06 & 0.16 & 1.00 & 0.30 & 1.00 & 0.44 & 1.00 & 0.44 & 1.00 & 0.47 & 0.0292\\
\hline
\end{tabular*}
\end{table}

We also included a standard linear regression (Lik) and Gervini and
Yohai's estimates (Gervini and Yohai \cite{Gervini2002}) based on a Huberized
estimate with an adaptively-chosen threshold (G--Y). Table~\ref
{linregtable} reports the means and standard deviations of the
parameters in this model calculated from 5000 simulations before
bootstrap methods are applied. We present computation times here as
well; bootstrapping results in multiplying these times by 100 for all
estimators.

As can be observed from these results, the use of multivariate density
estimation creates significant biases, particularly in $\beta_2$ and
$\beta_3$. This is mitigated in the centered density estimates,
although not for the homoscedastic estimators. We speculate that this
is because the conditional density estimate can correct for biasses
from the Nadaraya--Watson estimator which the homoscedastic restriction
does not allow for. The marginal methods perform considerably better
and achieve similar performance to those of Gervini and
Yohai \cite{Gervini2002}.
We also observe that Hellinger distance estimators have large variances
in some cases, mostly due to occasional outlying parameter estimates.
By contrast, negative exponential disparity estimators were much more stable.

In addition to the simulations above, for each simulated data set we
performed 100 bootstrap replicates as described in Section~\ref
{secbootstrap} and used this to both provide a bias correction and
confidence intervals. The resulting point estimates and coverage
probabilities are reported in Table~\ref{linregboottab}. Here we see
that much of the bias has been removed for all estimators except for
the homoscedastic models. Coverage probabilities are at least as close
to nominal values as minimum squared error estimators.

%
\begin{table}[b]
\tabcolsep=0pt
\caption{Statistical properties (estimate, standard deviation (sd) and
coverage (cov)) of inference following bootstrap bias correction and
using bootstrap confidence intervals. Labels for estimators are the
same as in Table~\protect\ref{linregtable}} \label{linregboottab}
\begin{tabular*}{\tablewidth}{@{\extracolsep{\fill}}@{}lllllllllllll@{}}
\hline
& $\beta_0^c$ & sd & cov & $\beta_1^c$ & sd & cov & $\beta_2^c$ & sd & cov & $\beta_3^c$ & sd & cov \\
\hline
Lik & 1.00 & 0.28 & 0.92 & 1.00 & 0.4 & 0.92 & 1.00 & 0.4 & 0.91 & 0.99 & 0.41 & 0.92 \\
Hub & 1.00 & 0.29 & 0.91 & 1.00 & 0.41 & 0.91 & 1.00 & 0.41 & 0.91 & 0.99 & 0.43 & 0.91 \\
G--Y & 1.01 & 0.31 & 0.91 & 1.00 & 0.43 & 0.92 & 1.00 & 0.43 & 0.91 & 0.99 & 0.45 & 0.92 \\
HD.c & 1.00 & 0.31 & 0.95 & 0.99 & 0.43 & 0.93 & 0.99 & 0.43 & 0.93 & 1.00 & 0.46 & 0.95 \\
NED.c & 1.00 & 0.31 & 0.96 & 0.99 & 0.42 & 0.94 & 0.99 & 0.43 & 0.94 & 1.00 & 0.46 & 0.95 \\
HD & 0.99 & 0.42 & 0.98 & 0.98 & 0.59 & 0.95 & 0.97 & 0.48 & 0.94 & 1.02 & 0.62 & 0.98 \\
NED & 0.99 & 0.3 & 0.98 & 0.97 & 0.4 & 0.96 & 0.97 & 0.4 & 0.96 & 1.02 & 0.44 & 0.98 \\
HD.h & 0.99 & 0.58 & 0.81 & 0.97 & 0.38 & 0.86 & 0.98 & 0.39 & 0.84 & 1.02 & 0.76 & 0.81 \\
NED.h & 1.00 & 0.36 & 0.86 & 0.97 & 0.39 & 0.86 & 0.98 & 0.4 & 0.86 & 1.01 & 0.52 & 0.87 \\
HD.m & 1.00 & 0.31 & 0.9 & 0.99 & 0.44 & 0.95 & 1.00 & 0.46 & 0.95 & 1.00 & 0.47 & 0.95 \\
NED.m & 1.00 & 0.32 & 0.9 & 0.99 & 0.46 & 0.94 & 1.00 & 0.47 & 0.94 & 1.00 & 0.48 & 0.95 \\
\hline
\end{tabular*}
\end{table}

\begin{figure}

\includegraphics{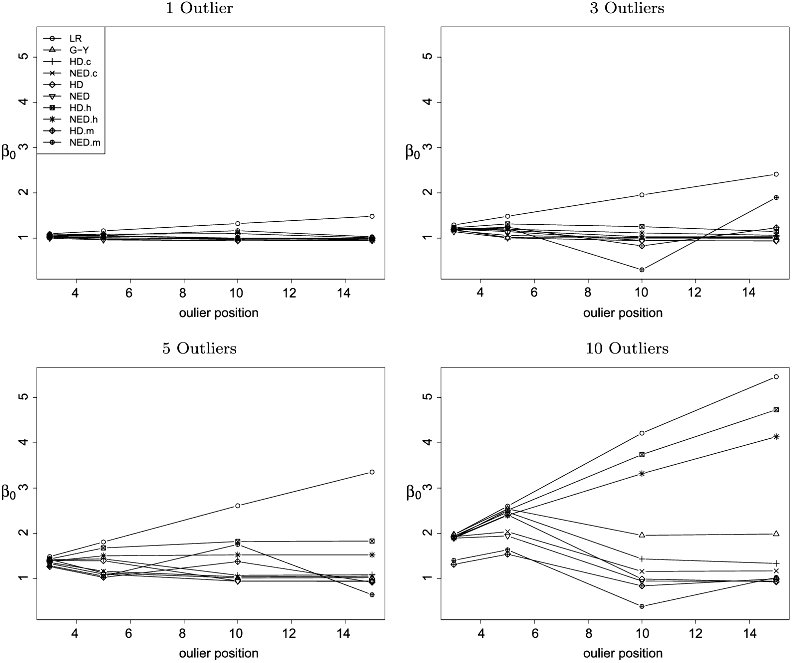}

\caption{Mean estimates $\hat{\beta}_0$ with different levels of
contamination uniformly distributed over
covariate values. Each line corresponds do a different estimation
method as given in the key.} \label{figO1beta0}
\end{figure}

To examine results when the data are contaminated, we plot the mean
estimate for $\beta_0$ under the contamination model~1 in Figure~\ref{figO1beta0} as the position of the
contamination increases; this mimics the bias plots of Lindsay \cite{Lindsay94}, Figures~1~and~2. We have reported plots at each level of
the number of contaminated observations. Here, we observe that the
least squares estimator is strongly affected although most robust
estimators are not. At 10 (30\%) contaminated observations, the
Gervini--Yohai estimator exhibits greater distortion of all except the
homoscedastic and maximum likelihood estimators, although it remains
robust and the tendency to ignore large outliers is evident. We
speculate that the breakdown in the homoscedastic methods is because
the underlying Nadaraya--Watson estimator is locally influenced strongly
by these values and the homoscedastic restriction does not allow it to
compensate for this. Estimates for the variance $\sigma$ were
similarly affected but the other regression parameters were not
influenced by outliers since they were uniformly distributed over the
range of covariates. A complete set of graphs is given in Figure~1 in Supplemental Appendix~A
(Hooker \cite{Hooker14}).

%

%
\begin{figure}

\includegraphics{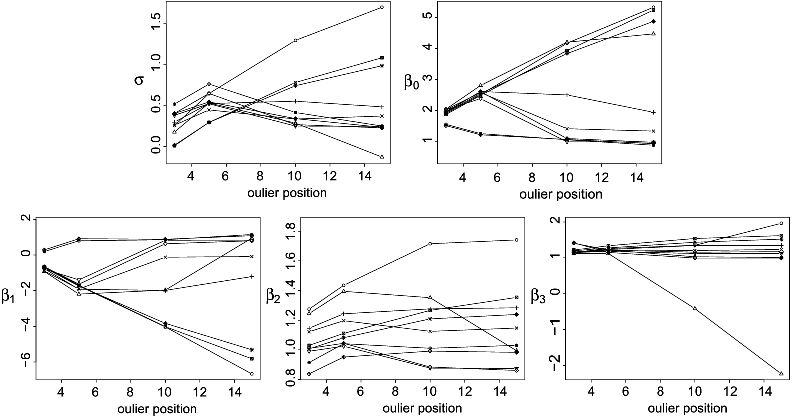}

\caption{Mean parameter estimates with 10 outliers with values $x_1$
close to $-$0.5.}\vspace*{-6pt} \label{figO2pars}
\end{figure}

By contrast, under contamination model~\ref{contamination2}, all
least-squares parameter estimates were affected by outliers. We have
plotted the average estimates for each parameter for 10 outliers in
Figure~\ref{figO2pars} using the same key as in Figure~\ref
{figO1beta0}. Here we observe that most estimators were robust,
although the Gervini--Yohai as well as the homoscedastic models were
affected. Investigating this more closely, at this level of
contamination, sampling distribution the Gervini--Yohai estimator
appears multi-modal which we speculate is associated with the adaptive
choice of the Huber threshold failing to reject some of the outliers.
It should be noted that this behavior was not evident at smaller
contamination percentages. Examining Figure~2 in
Supplemental Appendix~A (Hooker \cite{Hooker14}), we observe that this breakdown in robustness occurs most
dramatically only at 10 outliers, although the homoscedastic estimators
(but not Gervini--Yohai) show some evidence for this at 5 outliers as well.

\subsection{Logistic regression}

For logistic regression there is no option to center the response
before producing a conditional density estimate. We therefore examine
only the logistic regression (Lik), Hellinger distance (HD) and
negative exponential disparity (NED) estimators. Because logistic
regression estimates are less stable than linear regression, we used
121 points generated as described above. We also note that Monte Carlo
estimates are not required to evaluate the disparity in this case since
it is defined as a sum over a discrete set of points. Simulation
results are reported in Table~\ref{logistictable}.

%
\begin{table}[b]
\tabcolsep=0pt
\caption{Simulation results for logistic regression using maximum
likelihood (\texttt{LR}), Hellinger distance (\texttt{HD}) and negative
exponential disparity (NED) estimates}\label{logistictable}
\begin{tabular*}{\tablewidth}{@{\extracolsep{\fill}}@{}llllllllll@{}}
\hline
& $\beta_0$ & sd & $\beta_1$ & sd & $\beta_2$ & sd & $\beta_3$ & sd & Time \\
\hline
LR & \phantom{$-$}0.00 & 0.29 & 0.53 & 0.42 & 0.52 & 0.42 & 0.52 & 0.44 & 0.01 \\
HD & $-$0.01 & 0.33 & 0.57 & 0.44 & 0.56 & 0.44 & 0.58 & 0.49 & 0.01 \\
NED & $-$0.01 & 0.29 & 0.51 & 0.39 & 0.5 & 0.39 & 0.54 & 0.44 & 0.01 \\
\hline
\end{tabular*}
\end{table}

There is again a noticeable bias in these estimates and we employed the
bootstrapping methods outlined above both to remove the bias in the
estimates and to estimate confidence intervals. For each data set, we
simulated 100 bootstrap samples and used these to estimate the bias and
standard deviation of the estimators. In addition to removing bias, we
examined the coverage of a parametric bootstrap interval based on the
bias corrected estimate plus or minus 1.96 the bootstrap standard
deviation. The results of these experiments are reported in Table~\ref
{logisticboottable} where we observe that the bias has effectively
been removed, the standard deviations between the corrected estimators
are very similar between the disparity methods and standard logistic
regression estimates and we retain appropriate coverage levels.

%
\begin{table}
\tabcolsep=0pt
\caption{Simulation results for logistic regression following a
bootstrap to correct for bias and construct confidence intervals using
maximum likelihood \texttt{LR}, Hellinger distance \texttt{HD} and
negative exponential disparity (NED) estimates with mean estimate,
standard deviation across simulations (sd) and coverage of bootstrap
confidence intervals (cov)}\label{logisticboottable}
\begin{tabular*}{\tablewidth}{@{\extracolsep{\fill}}@{}lllllllllllll@{}}
\hline
& $\beta_0^c$ & sd & cov & $\beta_1^c$ & sd & cov & $\beta_2^c$ & sd & cov & $\beta_3^c$ & sd & cov \\
\hline
LR & $-$0.01 & 0.28 & 0.97 & 0.51 & 0.4 & 0.97 & 0.49 & 0.4 & 0.97 & 0.5 & 0.42 & 0.97 \\
HD & \phantom{$-$}0.00 & 0.31 & 0.96 & 0.55 & 0.43 & 0.95 & 0.53 & 0.44 & 0.94 & 0.52 & 0.46 & 0.96 \\
NED & $-$0.01 & 0.28 & 0.95 & 0.5 & 0.4 & 0.94 & 0.49 & 0.4 & 0.94 & 0.5 & 0.42 & 0.95 \\
\hline
\end{tabular*}
\end{table}

\begin{sidewaysfigure}

\includegraphics{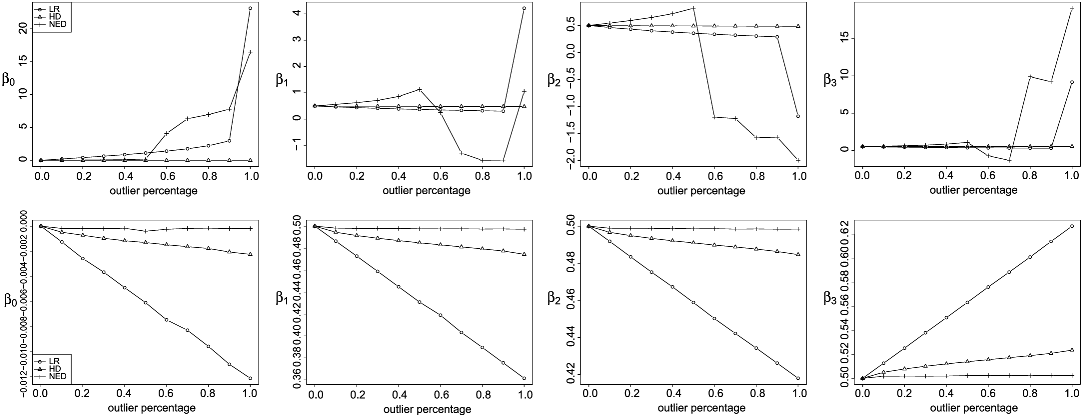}

\caption{Mean estimates of parameters in a logistic regression as the
outlier percentage increases. Top row: outliers occur uniformly over $X$.
Bottom: outliers at a single value of $X$.}
\label{figlogitRobust}
\end{sidewaysfigure}

The robustness of these estimates for binomial data from 8 trials at
each $X_i$ is examined in Figure~\ref{figlogitRobust}. Here we
observe that adding outliers at a single point generate classical
robust behavior -- the maximum likelihood estimate (calculating by
minimizing the Kullback--Leibler divergence) is highly non-robust while
Hellinger distance and negative exponential disparity are largely
unchanged. When outliers are added uniformly, we observe more
distortion of our estimates, particularly NED. This is both due to the
large over-all amount of contamination (at all points rather than just
one) and because we cannot achieve exact orthogonality between the
generating and contaminating distributions. At $\alpha= 0.5$, there
is, as expected, a significant change and both NED and HD exhibit
increased distortion.

%

\section{Real data} \label{secdata}

We demonstrate these methods with the analysis of the phosphorus
content data in \cite{RousseeuwLeroy05} in which plant
phosphorus in corn is related to organic and non-organic phosphorus in
the soil in which it is grown. In these data there is a distinct
outlier that significantly affects least squares estimates. However
robust procedures all produced estimates of approximately the same
magnitude. We also conducted a bootstrap analysis, as described in
Section~\ref{secsimulations} based on 100 bootstrap samples. The
results of these are reported in Table~\ref{phosphorustable}.

%
\begin{table}
\tabcolsep=0pt
\caption{Results on phosphorous data. Estimates with superscripts
($\beta^c$) incorporate a bootstrap bias correction, standard
deviations are also estimated via a bootstrap}\label{phosphorustable}
\begin{tabular*}{\tablewidth}{@{\extracolsep{\fill}}@{}lllllllllllll@{}}
\hline
&$\log\sigma$ & $\log\sigma^c$ & sd & $\beta_0$ & $\beta_0^c$ & sd & $\beta_1$ & $\beta_1^c$ & sd & $\beta_2$ & $\beta_2^c$ & sd \\
\hline
LR & 20.68 & 17.01 & 7.89 & 56.25 & 35.98 & 19.52 & 1.79 & 1.8 & 0.65 & \phantom{$-$}0.09 & \phantom{$-$}0.08 & 0.5 \\
Hub & \phantom{0}2.14 & \phantom{0}2.14 & 0.57 & 59.08 & 59.99 & 10.87 & 1.36 & 1.4 & 0.39 & \phantom{$-$}0.09 & \phantom{$-$}0.06 & 0.28 \\
G--Y & \phantom{0}2.51 & \phantom{0}2.8 & 0.38 & 66.47 & 63.02 & \phantom{0}8.86 & 1.29 & 1.28 & 0.33 & $-$0.11 & $-$0.05 & 0.23 \\
HD.c & \phantom{0}2.26 & \phantom{0}2.23 & 0.12 & 54.27 & 53.84 & \phantom{0}5.39 & 1.3 & 1.22 & 0.33 & \phantom{$-$}0.24 & \phantom{$-$}0.25 & 0.12 \\
NED.c & \phantom{0}2.16 & \phantom{0}2.14 & 0.12 & 53.19 & 53.08 & \phantom{0}6.78 & 1.23 & 1.15 & 0.32 & \phantom{$-$}0.27 & \phantom{$-$}0.27 & 0.15 \\
HD & \phantom{0}2.44 & \phantom{0}2.39 & 0.13 & 61.39 & 59.57 & 10.95 & 1.01 & 1.12 & 0.27 & \phantom{$-$}0.09 & \phantom{$-$}0.1 & 0.21 \\
NED & \phantom{0}2.4 & \phantom{0}2.4 & 0.16 & 56.78 & 52.45 & 14.08 & 1.03 & 1.15 & 0.3 & \phantom{$-$}0.19 & \phantom{$-$}0.25 & 0.26 \\
HD.h & \phantom{0}2.42 & \phantom{0}2.45 & 0.2 & 50.8 & 44.02 & 10.33 & 1.47 & 1.53 & 0.32 & \phantom{$-$}0.2 & \phantom{$-$}0.22 & 0.25 \\
NED.h & \phantom{0}2.33 & \phantom{0}2.32 & 0.18 & 52.77 & 49.08 & 10.29 & 1.35 & 1.31 & 0.3 & \phantom{$-$}0.21 & \phantom{$-$}0.24 & 0.26 \\
HD.m & \phantom{0}2.35 & \phantom{0}2.17 & 0.33 & 74.71 & 70.99 & 13.69 & 1.58 & 1.08 & 1.09 & $-$0.42 & $-$0.22 & 0.45 \\
NED.m & \phantom{0}2.36 & \phantom{0}2.28 & 0.32 & 60.33 & 57.2 & 11.46 & 1.21 & 1.08 & 0.71 & \phantom{$-$}0.1 & \phantom{$-$}0.22 & 0.35 \\
\hline
\end{tabular*}
\end{table}

\section{Discussion}

Conditionally specified models make up a large subset of the models
most commonly used in applied statistics, including regression,
generalized linear models and tabular data. In this paper, we
investigate the use of disparity methods to perform parameter
estimation across a range of such models. Our treatment is general in
covering multivariate response and covariate variables and allowing for
both discrete and continuous elements of each and almost any
probabilistic relationship between them. We have also investigated the
use of centering continuous responses by a Nadaraya--Watson estimator
based on a subset of the covariates and presented a complete theory
covering all ways to divide covariates into centering and conditioning
variables. Along the way we have established uniform $L_1$ convergence
results for a class of non-parametric conditional density estimates as
well as the consistency and a central limit theorem for disparity-based
models. These theoretical results highlight the consequences of
different choices of density estimate and disparity when the model is
incorrectly specified and demonstrate the limitations of centering
densities within this methodology unless the same covariates are used
within both the centering estimate and to condition. We have also
established a bootstrap bias correction and inference methodology that
has sound theoretical backing.

There are many direction for future study, starting from these methods.
As is the case for disparity estimators for multivariate data, the use
of conditional kernel densities results in a bias in parameter
estimates that cannot be ignored in our central limit theorem, except
in special cases. Empirically, our bootstrap methods reduce this bias,
but more sophisticated alternatives are possible. We have not
investigated using alternatives to Nadaraya--Watson estimators, but
conjecture that doing so may also reduce bias. In a linear regression
model, for example, the use of a local linear smoother should
completely remove the bias from $\hatmnh$ when the model is true. More
generally, centering based on a localized version of the assumed
parametric model may be helpful. An alternative method of removing the
bias follows the marginal approaches explored in Hooker and Vidyashankar \cite
{HookerVidyashankar13}. In this approach, the non-parametric density
estimate becomes dependent on a parametric transformation of the data
that is chosen in such a way that at the true parameters the
transformed data have independent dimensions. This would allow the use
of univariate density estimates, thereby removing the asymptotic bias.

In our examples, we have employed cross-validated log likelihood to
choose bandwidths and the robustness of this choice has not been
investigated. We speculate that a form of weighted cross-validation may
produce more robust bandwidth selection. We have also focussed solely
on kernel-based methods; little is known about the use of alternative
density estimates in disparity measures, although see Wu and
Hooker \cite
{WuHooker13} for an exploration of non-parametric Bayesian methods
combined with disparities.

Empirically, our methods perform very well in both the precision and
robustness of our estimators. Within our experiments, NED generally
improved upon HD methods; we speculate this is due to Hellinger
distance's sensitivity to inliers (see Lindsay \cite{Lindsay94}) and
hence added variability if the non-parametric estimate is sometimes
multi-modal. Moreover, in distinction to alternatives, our methods
provide a generic means of obtaining both robustness and efficiency
across a very wide range of applicable regression models.

The need for kernel density estimates for responses and covariates at
each level of the combined categorical variables limits the set of
situations in which our estimates are feasible at realistic sample
sizes. They are nonetheless relevant for non-trivial practical problems
in data analysis; the marginal approaches in Hooker and Vidyashankar \cite
{HookerVidyashankar13} also represent a means of approaching
higher-dimensional covariate spaces. These results open the way for the
application of minimum disparity estimates to a wide range of
real-world data analysis problems.

\begin{supplement}
\stitle{Proofs and simulations for consistency, efficiency and
robustness of conditional disparity methods}
\slink[doi]{10.3150/14-BEJ678SUPP} 
\sdatatype{.pdf}
\sfilename{BEJ678\_supp.pdf}
\sdescription{We provide additional supporting simulations of the
efficiency and robustness of the conditional disparity methods along
with proofs of the results stated above.}
\end{supplement}

\section*{Acknowledgements}
Research supported in part by NSF Grants DEB-0813743, CMG-0934735
and DMS-1053252. The author thanks Anand Vidyashankar for many helpful
discussions.



%

\printhistory
\end{document}